\documentclass[12pt]{article}
\usepackage{amsmath}
\usepackage{amssymb}
\usepackage{a4}
\usepackage{amscd}

\newcommand{\ssp}{\discretionary{}{}{\,}}

\newcounter{cdef}[section]
\newcounter{cthm}[section]
\newcounter{csatz}[section]
\newcounter{clem}[section]
\newcounter{cfolg}[section]
\newcounter{cbsp}[section]

\def\thecdef{\thesection.\arabic{cdef}}
\def\thecthm{\thesection.\arabic{cthm}}
\def\thecsatz{\thesection.\arabic{csatz}}
\def\theclem{\thesection.\arabic{clem}}
\def\thecfolg{\thesection.\arabic{cfolg}}
\def\thecbsp{\thesection.\arabic{cbsp}}

\newenvironment{bew}[1][{}]{\par\addvspace{.5\baselineskip}\par
\textbf{\bewname{#1. }}}{\ \hspace*{\fill}\rule{1ex}{1ex}
\par\addvspace{.5\baselineskip}\par}

\newenvironment{defin}{\par\addvspace{.5\baselineskip}\par
\refstepcounter{cdef} \textbf{\definname~\thecdef.} }{
\par\addvspace{.5\baselineskip}\par}

\newenvironment{thm}{\par\addvspace{.5\baselineskip}\par
\refstepcounter{cthm} \textbf{\thmname~\thecthm.} \em }{\em
\par\addvspace{.5\baselineskip}\par}

\newenvironment{satz}{\par\addvspace{.5\baselineskip}\par
\refstepcounter{csatz} \textbf{\satzname~\thecsatz.} \em }{\em
\par\addvspace{.5\baselineskip}\par}

\newenvironment{lemma}{\par\addvspace{.5\baselineskip}\par
\refstepcounter{clem} \textbf{\lemmaname~\theclem.} \em }{\em
\par\addvspace{.5\baselineskip}\par}

\newenvironment{folg}{\par\addvspace{.5\baselineskip}\par
\refstepcounter{cfolg} \textbf{\folgname~\thecfolg.} \em }{\em
\par\addvspace{.5\baselineskip}\par}

\newenvironment{bem}{\par\addvspace{.5\baselineskip}\par
\textbf{\bemname.} }{\ \hspace*{\fill}\rule{1ex}{1ex}
\par\addvspace{.5\baselineskip}\par}

\newenvironment{bems}{\par\addvspace{.5\baselineskip}\par
\textbf{\bemsname.} }{\ \hspace*{\fill}\rule{1ex}{1ex}
\par\addvspace{.5\baselineskip}\par}

\newcommand{\Rda}{\hat{R}}

\newcommand{\Rdam}{\hat{R}^{-1}}

\newcommand{\linv}[1]{(#1)_{\mathrm{l}}}
\newcommand{\rinv}[1]{(#1)_{\mathrm{r}}}

\newcommand{\cont}{^{\mathrm{c}}}
\newcommand{\contr}{\mathrm{c}}
\newcommand{\dif}{\mathrm{d}}
\newcommand{\Mor}{\mathrm{Mor}}
\newcommand{\id}{\mathrm{id}}

\def\bewname{Proof}
\def\definname{Definition}
\def\thmname{Theorem}
\def\satzname{Proposition}
\def\lemmaname{Lemma}
\def\folgname{Corollary}
\def\bemname{Remark}
\def\bemsname{Remarks}
\def\bspname{Example}

\newcommand{\btxandshort}[1]{and}
\newcommand{\btxpagesshort}[1]{pp.}
\newcommand{\Btxinshort}[1]{In}
\newcommand{\btxphdthesis}[1]{phd-thesis}

\newcommand{\mcl}{\mathcal}
\newcommand{\mf}{\mathfrak}
\newcommand{\mbb}{\mathbb}

\newcommand{\A}{\mcl{A}}
\newcommand{\otA}{\otimes _\A}
\newcommand{\comp}{\mathbb{C}}
\renewcommand{\linv}[1]{(#1)_\mathrm{l}}
\renewcommand{\rinv}[1]{(#1)_\mathrm{r}}
\newcommand{\binv}[1]{(#1)_\mathrm{lr}}
\newcommand{\Ga}[1]{\Gamma _{#1}}
\newcommand{\Gp}{\Ga{+}}
\newcommand{\Gm}{\Ga{-}}
\newcommand{\Gaw}[2][{}]{\Gamma ^{\wedge #1}_{#2}}
\newcommand{\Gat}[2][{}]{\Gamma ^{\otimes #1}_{#2}}
\newcommand{\mt}{,}
\newcommand{\tet}{\otA }
\newcommand{\ctr}[3][{}]{\langle {#2}\mt {#3}\rangle _{#1}}
\newcommand{\ctrp}[2]{\ctr[+]{#1}{#2}}
\newcommand{\ctrm}[2]{\ctr[-]{#1}{#2}}
\newcommand{\ctrpm}[2]{\ctr[\pm ]{#1}{#2}}
\newcommand{\ctrmp}[2]{\ctr[\mp ]{#1}{#2}}
\newcommand{\paar}[2]{\mathrm{B}({#1},{#2})}
\newcommand{\Hod}[2]{*_{#1}^{#2}}
\newcommand{\HoL}[1]{\Hod{\mathrm{L}}{#1}}
\newcommand{\HoR}[1]{\Hod{\mathrm{R}}{#1}}
\newcommand{\HoLp}{\HoL{+}}
\newcommand{\HoLm}{\HoL{-}}
\newcommand{\HoRp}{\HoR{+}}
\newcommand{\HoRm}{\HoR{-}}
\newcommand{\HoLpm}{\HoL{\pm }}
\newcommand{\HoLmp}{\HoL{\mp }}
\newcommand{\HoRpm}{\HoR{\pm }}
\newcommand{\HoRmp}{\HoR{\mp }}
\newcommand{\kod}[2][{}]{\partial _{#1}^{#2}}
\newcommand{\kodL}[1]{\kod[\mathrm{L}]{#1}}
\newcommand{\kodR}[1]{\kod[\mathrm{R}]{#1}}
\newcommand{\kodLpm}{\kodL{\pm }}

\newcommand{\kodRpm}{\kodR{\pm }}
\newcommand{\kodpm}[1]{\kod[#1]{\pm }}
\newcommand{\Lap}[2][\tau ]{\mathbf{\Delta }^{#2}_{#1}}
\newcommand{\Lappm}[1][\tau ]{\Lap[#1]{\pm }}

\newcommand{\br}{\mathbf{r}}
\newcommand{\ract}{\triangleleft}
\newcommand{\Lto}[2]{{[{#1}\leftarrow {#2}]}}
\newcommand{\Rto}[2]{{[{#1}\rightarrow {#2}]}}
\newcommand{\RLto}[2]{{[{#1}\rightleftarrows {#2}]}}
\newcommand{\SLqN}{\mathrm{SL}_q(N)}
\newcommand{\OSLqN}{\mcl{O}(\mathrm{SL}_q(N))}
\newcommand{\OSpqN}{\mcl{O}(\mathrm{Sp}_q(N))}
\newcommand{\OOqN}{\mcl{O}(\mathrm{O}_q(N))}
\newcommand{\komult}{\varDelta }
\newcommand{\kowl}{\varDelta _\mathrm{L}}
\newcommand{\kowr}{\varDelta _\mathrm{R}}

\title{Hodge and Laplace-Beltrami Operators for
Bicovariant Differential Calculi on Quantum Groups}
\author{Istv\'an Heckenberger
\thanks{e-mail: heckenbe@mathematik.uni-leipzig.de}\\
\textit{Universit{\"a}t Leipzig, Mathematisches Institut}}
\date{}

\begin{document}
\maketitle

\begin{abstract}
For bicovariant differential calculi on quantum matrix groups a
generalisation of classical notions such as metric tensor, Hodge operator,
codifferential and Laplace-Beltrami
operator for arbitrary $k$-forms is given. Under some technical assumptions
it is proved that Woronowicz' external algebra of left-invariant differential
forms either contains a unique form of maximal degree or it is infinite
dimensional. Using
Jucys-Murphy elements of the Hecke algebra the eigenvalues of
the Laplace-Beltrami operator for the Hopf algebra $\OSLqN$ are computed.
\end{abstract}

\section{Introduction}
\label{sec-einfuehrung}

About ten years ago S.\,L.~Woronowicz introduced the concept of bicovariant
differential calculus on arbitrary Hopf algebras and developed a general
theory of such calculi \cite{a-Woro2}. One of the most interesting parts
of this theory is his definition of external algebras and higher order
calculi by using a braiding map instead of the flip operator in the
corresponding classical constructions. The higher order differential
calculus defined in this manner becomes then an $\mbb{N}_0$-graded
differential super Hopf algebra (\cite{a-Brz1}; see \cite{b-KS} for a complete
proof). However, applying Woronowicz' construction of higher order calculi
to quantum matrix groups leads to a number of difficulties and
phenomena that do not occur in the classical (commutative) case.
Firstly, the vector space $\linv{\Gaw{}}$ of left-invariant differential
forms endowed with the canonical (wedge) product does not form a Grassmann
algebra in general.
Secondly, it may happen that the dimensions of the spaces $\linv{\Gaw[k]{}}$
of left-invariant $k$-forms do not vanish as $k\to +\infty $
(see \cite{a-HSchu2}). For the irreducible $N^2$-dimensional bicovariant
first order differential calculi on the coordinate Hopf algebra $\OSLqN $
of the quantum group $\SLqN $, $N\geq 2$,
a detailed description of the higher order differential calculi $\Gaw{}$
was given by A.~Sch{\"u}ler \cite{a-Schueler1}. In this important case
it is proved in \cite{a-Schueler1} that for transcendental values of the
parameter $q$ the dimension of the vector space of left-invariant k-forms
is $\binom{N^2}{k} $ just as in the classical situation.

In ``ordinary" differential geometry the Laplace-Beltrami operator
$\Lap[]{}$ acting on differential forms plays a central role.
In its construction a metric tensor, the Hodge star and the codifferential
operators are essentially used.
The aim of this paper is to give a definition of invariant
Laplace-Beltrami operators $\Lap[]{}$ for inner bicovariant differential
calculi on arbitrary Hopf algebras. It will be a generalisation of the
classical concept and works also in the case
when the higher order calculus is infinite dimensional. The existence of
$\Lap[]{}$ is shown for coquasitriangular Hopf algebras and irreducible
differential calculi defined by generalised $l$-functionals.
As tools we use $\sigma $-metrics (a generalisation of the concept of
a metric tensor in the commutative case), Hodge star operators
(in a special case) and codifferentials.

In Section \ref{sec-metriken} we introduce $\sigma $-metrics for a pair
of bicovariant bimodules. In Section \ref{sec-examples} we give examples
for these structures.
In Section \ref{sec-kontraktionen} further basic notions
like contractions with forms (see also \cite{a-DurdOzi}) and
$\sigma $-metrics on higher order forms of Woronowicz' external algebra
are introduced and a number of useful properties of these mappings
are developed.
Section \ref{sec-Hodge} is concerned with Hodge operators and codifferential
operators. For their definitions we require two assumptions. The first one
is that the Hopf algebra is ``connected"
(i.\,e.\ it has only one one-dimensional corepresentation), and the second
assumption is satisfied (for instance) if the left-invariant part of
the external algebra is finite dimensional.
In Theorem \ref{t-eindhf} it is proved that if there is a left-covariant
$\sigma$-metric on the external algebra then there exists a unique
(up to a complex multiple) left-invariant differential form of maximal
degree. For the proof of Theorem \ref{t-eindhf} (and its Corollary
\ref{f-eindhf})
we don't need the assumption that the Hopf algebra is ``connected".
Further we define Hodge star and codifferential operators and prove some
of their properties. One of the formulas for the codifferential operator is
independent of the Hodge star and will be taken as a definition in
the next section. In Section \ref{sec-Laplace} the invariant
Laplace-Beltrami operator is defined and a number of results on this operator
are derived. Among others, it is shown (Theorem \ref{t-hdifkod}) that there
is a duality between the differential and codifferential as in the
classical case.
In Section \ref{sec-eigenvalues} the eigenvalues of the Laplace-Beltrami
operator for the quantum group $\SLqN$, $N\geq 2$ are determined.

In this paper we shall use the convention to sum over repeated indices
belonging to different terms.
Throughout, $\A $ denotes a Hopf algebra over the complex field
with comultiplication $\komult $ and invertible antipode $S$.
The symbol $\otA $ means the algebraic tensor product over the Hopf algebra
$\A $, while $\kowl $ and $\kowr $ denote left and right coactions on a
bicovariant $\A $-bimodule, respectively. If $u$ and $v$ are
corepresentations of $\A $,
then we write $\Mor(u,v)$ for the set of intertwiners of $u$ and $v$.
We set $\Mor(u)=\Mor(u,u)$. Throughout the paper we freely use basic facts
from the theory of bicovariant differential calculi (see \cite{a-Woro2}
or \cite[Chapter 14]{b-KS}.

I want to thank Prof.~Schm\"udgen for posing the problem and for motivating
discussions.

\section{$\sigma $-Metrics}
\label{sec-metriken}

Let $\A $ be an arbitrary Hopf algebra and let $\Gp $ and $\Gm $ be two
finite dimensional bicovariant $\A $-bimodules. Recall that
any bicovariant bimodule $\Ga{}$ is a free left and right $\A $-module and
there are bases of $\Ga{}$ consisting of left- and right-invariant
elements respectively.
In what follows we use the symbols
$\linv{\Ga{}},\rinv{\Ga{}}$ and $\binv{\Ga{}}$ to denote the vector spaces
of left-, right- and biinvariant (i.\,e.\ both left- and right-invariant)
elements in a bicovariant bimodule $\Ga{}$.
Further, there is a canonical braiding $\sigma :\Ga{\tau } \otA \Ga{\tau '}
\to \Ga{\tau '}\otA \Ga{\tau }$ (defined by Woronowicz \cite{a-Woro2}) for
each $\tau ,\tau '\in \{+,-\}$ which is an invertible homomorphism
of bicovariant bimodules. We shall write $\sigma ^+$ for $\sigma $
and $\sigma ^-$ for $\sigma ^{-1}$.

\begin{defin}\label{d-metrik}
A linear mapping $g:\Gp \otA \Gm+\Gm\otA \Gp \to \A $ is called a
\textit{$\sigma $-metric of the}
(not ordered) \textit{pair} $(\Gp ,\Gm )$ if it satisfies the following
conditions:
\begin{itemize}
\item
$g$ is a homomorphism of $\A $-bimodules,
\item
$g$ is nondegenerate, (i.\,e.\ for $\xi \in \Ga{\tau }$ both
$g(\xi \otA \xi ')=0$ for any $\xi '\in \Ga{-\tau }$ and
$g(\xi '\otA \xi )=0$ for any $\xi '\in \Ga{-\tau }$ imply $\xi =0$)
\item
$g\circ \sigma =g$ ($\sigma $-symmetry),
\item
the following diagrams commute ($\tau ,\tau '\in \{+,-\}$):
\begin{gather}\label{eq-g12sig23}
\begin{CD}
\Ga{\tau } \otA \Ga{\tau '} \otA \Ga{-\tau } @>\sigma ^\pm _{23}>>
\Ga{\tau } \otA \Ga{-\tau } \otA \Ga{\tau '}\\
@V\sigma ^\mp _{12}VV @VVg_{12}V\\
\Ga{\tau '} \otA \Ga{\tau } \otA \Ga{-\tau } @>g_{23}>>\Ga{\tau '}
\end{CD}
\end{gather}
\end{itemize}
The $\sigma $-metric of the pair $(\Gp ,\Gm )$ is said to be
\textit{left-covariant} resp.\ \textit{right-covariant} if
\begin{align}\label{eq-lkmet}
\komult \circ g&=(\id \otimes g)\kowl \quad\text{resp.\ }\\
\komult \circ g&=(g\otimes \id )\kowr
\end{align}
on $\Gp \otA \Gm +\Gm \otA \Gp $. We call it \textit{bicovariant}
if it is both left- and right-covariant.
\end{defin}
If no ambiguity can arise then we use the symbol `,' in order to separate
the two arguments of $g$. Recall that by definition we still have
$g(\xi a,\rho )=g(\xi ,a\rho )$ for any $a\in \A $,
$\xi \in \Ga{\tau }$ and $\rho \in \Ga{-\tau }$, $\tau \in \{+,-\}$.

If $g$ is a homomorphism of the $\A $-bimodules
$\Ga{\tau }\otA \Ga{-\tau }$ and $\A $, $\tau \in \{+,-\}$,
(e.\,g.\ if $g$ is a $\sigma $-metric of the pair $(\Gp ,\Gm )$)
then on the tensor product $\bigotimes _{m=1}^k\Ga{\tau _m}$,
$\tau _m\in \{+,-\}$ the equation
\begin{align}\label{eq-metvert}
g_{i,i+1}\circ g_{j,j+1}&=g_{j-2,j-1}\circ g_{i,i+1}
\quad \text{for $k>j>i+1$}
\end{align}
holds. One can check that the only conditions on the above map to be
well defined is $\tau _{i+1}=-\tau _i$ and $\tau _{j+1}=-\tau _j$.
The formulas
\begin{gather}
g_{i,i+1}\circ \sigma _{j,j+1}=\sigma _{j-2,j-1}\circ g_{i,i+1}
\quad \text{and}\quad
\sigma _{i,i+1}\circ g_{j,j+1}=g_{j,j+1}\circ \sigma _{i,i+1},
\quad \text{$j>i+1$,}
\end{gather}
should be clear as well.

Let now $g$ be a $\sigma $-metric of the pair $(\Gp ,\Gm )$. Then on the
tensor product $\Gat{\tau }\otA \Gat{-\tau }$,
$\tau \in \{+,-\}$ we define a map $\tilde{g}$
recursively by setting
\begin{equation}\label{eq-tildeg}
\begin{gathered}
\tilde{g}(\xi \mt a):=\xi a,\quad
\tilde{g}(a\mt \zeta):=a\zeta,\\
\tilde{g}(\xi \otA \xi _1\mt \zeta _1\otA \zeta ):=
\tilde{g}(\xi g(\xi _1\mt \zeta _1)\mt \zeta )
\end{gathered}
\end{equation}
for all $\xi \in \Gat{\tau }$, $\xi _1\in \Ga{\tau }$,
$\zeta \in \Gat{-\tau }$, $\zeta _1\in \Ga{-\tau }$ and $a\in \A $.

Since $g$ is a homomorphism of $\A $-bimodules,
the map $\tilde{g}$ is well defined and it is a homomorphism of
bimodules. Note that $\tilde{g}$ is left-, right- or
bicovariant if $g$ is. The next lemma is crucial in what follows.

\begin{lemma}\label{l-tgsig}
For a $\sigma $-metric $g$ of the pair $(\Gp ,\Gm )$ and arbitrary
integers $i,k,l$ such that $1\leq i<k,l$, we have
\begin{gather}
\tilde{g}\circ (\sigma ^\pm _{k-i,k-i+1}\mt \id^{\otimes l})
=\tilde{g}\circ (\id ^{\otimes k}\mt \sigma ^\pm _{i,i+1})
\end{gather}
on the bimodule $\Gat[k]{\tau }\otA \Gat[l]{-\tau }$.
\end{lemma}

\begin{bew}
Because of (\ref{eq-tildeg}) it suffices to show the assertion for
$i=1$ and $k=l=2$. But in this case we have $\tilde{g}=g_{12}\circ g_{23}$
and it suffices to apply the fourth condition on the $\sigma $-metric
$g$ (see (\ref{eq-g12sig23}) in Definition \ref{d-metrik}) twice. We obtain
\begin{align*}
\tilde{g}(\sigma ^\pm (\xi _1\otA \xi _2)\mt \zeta _1\otA \zeta _2)&=
g_{12}\circ g_{23}\circ \sigma ^\pm _{12}
(\xi _1\otA \xi _2\otA \zeta _1\otA \zeta _2)\\
&=g_{12}\circ g_{12}\circ \sigma ^\mp _{23}
(\xi _1\otA \xi _2\otA \zeta _1\otA \zeta _2)\\
&=g_{12}\circ g_{34}\circ \sigma ^\mp _{23}
(\xi _1\otA \xi _2\otA \zeta _1\otA \zeta _2)\\
&=g_{12}\circ g_{23}\circ \sigma ^\pm _{34}
(\xi _1\otA \xi _2\otA \zeta _1\otA \zeta _2)\\
&=\tilde{g}(\xi _1\otA \xi _2\mt \sigma ^\pm (\zeta _1\otA \zeta _2)),
\end{align*}
where the third equation follows from (\ref{eq-metvert}).
\end{bew}

Let $g$ be a homomorphism of the bicovariant bimodules
$\Gp \otA \Gm + \Gm \otA \Gp $ and $\A $.
The general theory of bicovariant bimodules assures that $g$ is
nondegenerate whenever the matrix of $g$ with respect to one
fixed basis of $\linv{\Gp }$ and one fixed basis of
$\linv{\Gm }$ is invertible.
Conversely, if $g$ is left-covariant (i.\,e.\ (\ref{eq-lkmet}) is
fulfilled) then the matrix $G$ of $g$ with respect to any basis of
$\linv{\Gp }$ and $\linv{\Gm }$ has complex entries and the
nondegeneracy of $g$ implies the invertibility of the matrix $G$.
In this case we easily conclude that the following assertions are
equivalent:
\begin{description}
\item
(i) $g$ is nondegenerate,
\item
(ii) the restriction of $g$ onto the subspace
$\linv{\Gp \otA \Gm }+\linv{\Gm \otA \Gp }$ is nondegenerate,
\item
(iii) the matrix $G$ of $g$ with respect to one (and then any) basis of
$\linv{\Gp }$ and $\linv{\Gm }$ is invertible.
\end{description}
Obviously, this holds for left-covariant $\sigma $-metrics as well.
In what follows most of the $\sigma $-metrics will be left-covariant.

\section{Examples}
\label{sec-examples}

Let $\A $ be a coquasitriangular Hopf algebra (see for example
\cite{b-KS}, Section 10.1) with universal $r$-form $\br $ and let
$u=(u^i_j)_{i,j=1,\ldots,d}$ be a corepresentation of $\A $.
Then $u\cont =((u\cont )^i_j)_{i,j=1,\ldots ,d}$, $(u\cont )^i_j=S(u^j_i)$
is the contragredient corepresentation of $u$ and $u$ and $u\cont $
determine two bicovariant $\A $-bimodules $\Gp $ and $\Gm $, respectively.
They are given
by fixing the bases $\{\omega _{ij}\,|\,i,j=1,\ldots,d\}$ and
$\{\theta _{ij}\,|\,i,j=1,\ldots,d\}$ of left-invariant forms of
$\Gp $ resp.\ $\Gm $ and defining the right coactions $\kowr $ and right
$\A $-actions $\xi \ract a=S(a_{(1)})\xi a_{(2)}$,
$\xi \in \linv{\Ga{\tau }}$, $\tau \in \{+,-\}$, $a\in \A $, by the formulas
\begin{gather}
\kowr (\omega _{ij})=\omega _{kl}\otimes (uu\cont )^{kl}_{ij},\quad
\kowr (\theta _{ij})=\theta _{kl}\otimes
(u^{\contr \contr }u\cont )^{kl}_{ij},
\end{gather}
\begin{align}
\omega _{ij}\ract a&=S(l^-{}^k_i)l^+{}^j_l(a)\omega _{kl}
=\br (u^k_i,a_{(1)})\br (a_{(2)},u^j_l)\omega _{kl},\\
\theta _{ij}\ract a&=l^+{}^i_kS(l^-{}^l_j)(a)\theta _{kl}
=\br (a_{(1)},S(u^k_i))\br (S(u^j_l),a_{(2)})\theta _{kl}.
\end{align}
Note that the 1-forms $\omega :=\sum _{i=1}^d\omega _{ii}\in \Gp $
and $\theta :=(f\circ S)(u^i_j)\theta _{ij}\in \Gm $, where
$f(a)=\br (a_{(1)},S(a_{(2)}))$, are biinvariant.

Assume for a moment that the corepresentations $u$ and $u\cont $
are equivalent ($u\cong u\cont $) and let $T=(T^i_j)_{i,j=1,\ldots,d}$
be an invertible
morphism $T\in \Mor(u,u\cont )$. Clearly we have $T^{-1}\in\Mor(u\cont,u)$.
Then the mapping
\begin{equation}\label{eq-mpisom}
\theta_{ij}\mapsto \br(u^k_r,u^s_l)(T^{-1})^r_jT^i_s\omega_{kl}
\end{equation}
extends uniquely to a homomorphism of the bicovariant bimodules
$\Gm$ and $\Gp$. Moreover, this mapping is invertible and its inverse
is given by
\begin{equation}\label{eq-pmisom}
\omega_{ij}\mapsto \br(u^r_i,S(u^j_s))T^l_r(T^{-1})^s_k\theta_{kl}.
\end{equation}
We also see easily that this isomorphism maps $\theta$ into $\omega$.

Let now $u$ be an arbitrary corepresentation and let
$F_1\in\Mor(u^{\contr\contr},u)$, $F_2\in\Mor(u,u^{\contr\contr})$ and
$G_1,G_2\in\Mor(u)$ be invertible morphisms. 
Then we define linear maps $g':\Gp\otA\Gm\to\A$ and
$g'':\Gm\otA\Gp\to\A$ by
\begin{equation}\label{eq-g1g2}
g'(a\omega_{ij}\otA\theta_{kl})=aF_1{}^j_kF_2{}^l_i\quad\text{ and }\quad
g''(a\theta_{ij}\otA\omega_{kl})=aG_1{}^j_kG_2{}^l_i.
\end{equation}

\begin{lemma}\label{l-g1g2}
The mappings $g':\Gp \otA \Gm \to \A $ and $g'':\Gm \otA \Gp \to \A $ are
homomorphisms of bicovariant bimodules. Moreover, as bilinear forms they are
nondegenerate.
\end{lemma}

\begin{bew}
Firstly let us show that
$g'(\omega_{ij}\otA\theta_{kl}a)=g'(\omega_{ij}\otA\theta_{kl})a$.
For this we compute
\begin{align*}
\allowdisplaybreaks
g'(\omega_{ij}\otA\theta_{kl}a)&=
g'(a_{(1)}\br(u^r_i,a_{(2)})\br(a_{(3)},u^j_s)\br(a_{(4)},S(u^p_k))
\br(S(u^l_n),a_{(5)})\omega_{rs}\otA\theta_{pn})\\
&=a_{(1)}\br(u^r_i,a_{(2)})\br(a_{(3)},u^j_s)\br(a_{(4)},S(u^p_k))
\br(S(u^l_n),a_{(5)})F_1{}^s_pF_2{}^n_r\\
&=a_{(1)}\br(F_2{}^n_ru^r_i,a_{(2)})\br(a_{(3)},S(u^p_k)u^j_sF_1{}^s_p)
\br(S(u^l_n),a_{(4)})\\
&=a_{(1)}\br(F_2{}^n_ru^r_i,a_{(2)})\br(a_{(3)},S(u^p_k)F_1{}^j_sS^2(u^s_p))
\br(S(u^l_n),a_{(4)})\\
&=F_1{}^j_sa_{(1)}\br(F_2{}^n_ru^r_i,a_{(2)})\br(a_{(3)},S(S(u^s_p)u^p_k))
\br(S(u^l_n),a_{(4)})\\
&=F_1{}^j_ka_{(1)}\br(F_2{}^n_ru^r_i,a_{(2)})\br(S(u^l_n),a_{(3)})\\
&=F_1{}^j_ka_{(1)}\br(S^2(u^n_r)F_2{}^r_iS(u^l_n),a_{(2)})\\
&=F_1{}^j_kF_2{}^r_ia_{(1)}\br(S(u^l_nS(u^n_r)),a_{(2)})\\
&=F_1{}^j_kF_2{}^l_ia=g'(\omega_{ij}\otA\theta_{kl})a.
\end{align*}
Secondly we prove the covariance of $g'$, that is
\begin{gather}\label{eq-kovmet}
(\id\otimes g')\kowl =\komult \circ g'\quad \text{and}\quad
(g'\otimes\id)\kowr =\komult \circ g'
\end{gather}
as a mapping from $\Gp\otA\Gm\to\A\otimes\A$.
Similarly to the proof of Lemma 2.1 in \cite{a-HeckSch1} one can show that
the equations (\ref{eq-kovmet}) are equivalent to
$g'(\omega _{ij}\otA \theta _{kl})\in \comp $
and $g'(\omega _{ij}\otA \theta _{kl})
(uu\cont u^{\contr \contr }u\cont )^{ijkl}_{mnrs}
=g'(\omega _{mn}\otA \theta _{rs})$.
The first one is trivial. For the second we compute
\begin{align*}
&g'(\omega_{ij}\otA\theta_{kl})(uu\cont u\cont{}\cont u\cont)^{ijkl}_{mnrs}
=F_1{}^j_k F_2{}^l_i u^i_m S(u^n_j) S^2(u^k_r) S(u^s_l)\\
&\qquad=F_2{}^l_i u^i_m S(u^n_j) u^j_k F_1{}^k_r S(u^s_l)
=F_1{}^n_r S^2(u^l_i) F_2{}^i_m S(u^s_l)\\
&\qquad=F_1{}^n_r F_2{}^i_m S(u^s_lS(u^l_i))
=F_1{}^n_r F_2{}^s_m=g'(\omega _{mn}\otA \theta _{rs}).
\end{align*}
Hence the assertion follows.

Thirdly we have to prove the nondegeneracy of $g'$.
We shall carry out the proof only for the second argument of $g'$.
Let $\rho $ be an arbitrary element of $\Gm $. Then there are elements
$a_{ij}\in \A $ such that $\rho =\theta _{ij}a_{ij}$. Assume that
$g'(\rho' \otA \rho )=0$ for all $\rho '\in \Gp $. Inserting
$\rho '=\omega _{kl}$, $k,l=1,\ldots ,d$ and using that $g'$ is a right
$\A $-linear mapping, we obtain $F_1{}^l_i F_2{}^j_k a_{ij}=0$
for all $k,l$ and the invertibility of $F_1$ and $F_2$ gives
$a_{ij}=0$. Hence $\rho =0$.
\end{bew}

Assume for a moment that the corepresentations $u$ and $u\cont $ are
equivalent and let us identify $\Gm $ and $\Gp $ via the
isomorphism (\ref{eq-mpisom}).

\begin{lemma}\label{l-pmmetrik}
Suppose that the corepresentations $u$ and $u\cont $ are equivalent.
Then the homomorphisms $g'$ and
$g'':\Gp \otA \Gp $ in Lemma \ref{l-g1g2} coincide if and only if
there is a nonzero complex number $c$ such that
\begin{equation}
F_1{}^i_j=c(T^{-1})^i_rG_2{}^s_rT^j_s\quad\text{and}\quad
F_2{}^i_j=c^{-1}(T^{-1})^r_iG_1{}^s_rT^s_j.
\end{equation}
\end{lemma}

\begin{bew}
Since $g'$ and $g''$ are homomorphisms of $\A $-bimodules
it suffices to prove the assertion on the vector space
$\linv{\Gp }\otimes \linv{\Gp }$.
Inserting (\ref{eq-pmisom}) into the definition of $g'$ and $g''$,
it follows that
$g'(\omega _{ij}\otA \omega _{kl})=
g''(\omega _{ij}\otA \omega _{kl})$
if and only if
\begin{gather*}
T^y_r(T^{-1})^s_x\br (u^r_k,S(u^l_s))F_1{}^j_xF_2{}^y_i=
\br (u^r_i,S(u^j_s))T^y_r(T^{-1})^s_xG_1{}^y_kG_2{}^l_x
\end{gather*}
for any $i,j,k,l$. For the right hand side we compute
\begin{align*}
&\br (T^y_ru^r_i,S(u^j_s(T^{-1})^s_x))G_1{}^y_kG_2{}^l_x=
\br (S(u^r_y)T^r_i,S((T^{-1})^j_sS(u^x_s)))G_1{}^y_kG_2{}^l_x\\
&\qquad = T^r_i(T^{-1})^j_s\br (u^r_yG_1{}^y_k,S(G_2{}^l_xu^x_s))=
T^r_i(T^{-1})^j_s\br (G_1{}^r_yu^y_k,S(u^l_xG_2{}^x_s))
\end{align*}
and hence the lemma is valid if and only if
\begin{gather*}
T^y_r(T^{-1})^s_x\br (u^r_k,S(u^l_s))F_1{}^j_xF_2{}^y_i=
T^r_i(T^{-1})^j_s\br (G_1{}^r_yu^y_k,S(u^l_xG_2{}^x_s)).
\end{gather*}
Multiplying this equation by $\br (u^k_z,u^t_l)(T^{-1})^z_nT^m_t$
we obtain the equivalent condition
\begin{gather*}
F_1{}^j_mF_2{}^n_i=\bigl((T^{-1})^j_sG_2{}^t_sT^m_t\bigr)
\bigl((T^{-1})^z_nG_1{}^r_tT^r_i\bigr)
\end{gather*}
for any $i,j,m,n$, from which the assertion follows.
\end{bew}

Now let $u$ be an arbitrary corepresentation of $\A $
and let $g$ be the homomorphism on $\Gp \otA \Gm +\Gm \otA \Gp $
given by $g'$ and $g''$.
To prove the third and fourth conditions of Definition \ref{d-metrik} for $g$
let us recall the following explicit formulas for the braiding $\sigma $ 
(see \cite{b-KS}, Section 13.1):
\begin{align*}
\sigma (\omega _{ij}\otA \omega _{kl})&=
\br (u^r_t,S(u^y_n))\br (u^t_i,u^m_x)\br (S(u^l_y),u^z_s)\br (u^x_k,u^j_z)
\omega _{mn}\otA \omega _{rs},\\
\sigma (\omega _{ij}\otA \theta _{kl})&=
\br (u^r_t,S(u^y_n))\br (u^t_i,S^2(u^m_x))\br (S(u^l_y),u^z_s)
\br (S^2(u^x_k),u^j_z)\theta _{mn}\otA \omega _{rs},\\
\sigma (\theta _{ij}\otA \omega _{kl})&=
\br (u^y_n,u^r_t)\br (u^m_x,S(u^t_i))\br (u^z_s,u^l_y)
\br (S(u^j_z),u^x_k)\omega _{mn}\otA \theta _{rs},\\
\sigma (\theta _{ij}\otA \theta _{kl})&=
\br (u^y_n,u^r_t)\br (S(u^m_x),u^t_i)\br (u^z_s,u^l_y)
\br (u^j_z,S(u^x_k))\theta _{mn}\otA \theta _{rs}.
\end{align*}
The inverses $\sigma^{-1}$ of these braidings take the form
\begin{align*}
\sigma ^{-1}(\omega _{ij}\otA \omega _{kl})&=
\br (u^r_t,S(u^y_n))\br (S(u^m_x),u^t_i)\br (u^z_s,u^l_y,)\br (u^x_k,u^j_z)
\omega _{mn}\otA \omega _{rs},\\
\sigma ^{-1}(\omega _{ij}\otA \theta _{kl})&=
\br (S^2(u^y_n),u^r_t)\br (u^t_i,S^2(u^m_x))\br (S(u^l_y),u^z_s)
\br (u^j_z,S(u^x_k))\theta _{mn}\otA \omega _{rs},\\
\sigma ^{-1}(\theta _{ij}\otA \omega _{kl})&=
\br (S(u^r_t),u^y_n)\br (u^m_x,S(u^t_i))\br (u^z_s,u^l_y)
\br (u^x_k,u^j_z)\omega _{mn}\otA \theta _{rs},\\
\sigma ^{-1}(\theta _{ij}\otA \theta _{kl})&=
\br (u^y_n,u^r_t)\br (u^t_i,u^m_x)\br (S(u^l_y),u^z_s)
\br (u^j_z,S(u^x_k))\theta _{mn}\otA \theta _{rs}.
\end{align*}

\begin{satz}\label{s-4bedmet}
Let $F_1\in\Mor(u^{\contr\contr},u)$, $F_2\in\Mor(u,u^{\contr\contr})$ and
$G_1,G_2\in\Mor(u)$ be arbitrary invertible morphisms. 
Then the bilinear map $g:\Gp \otA \Gm +\Gm \otA \Gp $ given by $g'$ and
$g''$ in Lemma \ref{l-g1g2} satisfies the fourth condition
\begin{gather*}
g_{12}\sigma _{23}^\pm (\xi _1\otA \xi _2\otA \xi _3)=
g_{23}\sigma _{12}^\mp (\xi _1\otA \xi _2\otA \xi _3),
\end{gather*}
$\xi _1\in \Ga{\tau }$, $\xi _2\in \Ga{\tau '}$, $\xi _3\in \Ga{-\tau }$,
$\tau ,\tau '\in\{+,-\}$, of Definition \ref{d-metrik}.
\end{satz}

\begin{bew}
Since $g$ is a homomorphism of $\A $-bimodules (see Lemma \ref{l-g1g2})
it suffices to prove the assertion on the vector spaces
$\linv{\Ga{\tau }}\otimes \linv{\Ga{\tau '}}\otimes \linv{\Ga{-\tau }}$,
$\tau ,\tau '\in \{+,-\}$. We have to consider four cases
which correspond to the possible values of $\tau $ and $\tau '$.
Since the proofs are very similar, we only show the assertion
$g_{12}\sigma ^{-1}_{23}=g_{23}\sigma _{12}$
for $\tau =+$ and $\tau '=-$.
We will only use the formula
$\br (S(a),S(b))=\br (a,b)$ for any $a,b\in \A $ and the properties of
$F_i$ and $G_i$, $i=1,2$.
\begin{align*}
&g_{23}\bigl(\sigma (\omega _{ij}\otA \theta _{kl})\otA
\theta _{ab}\bigr)\\
&\quad=\br (u^r_t,S(u^y_n))\br (u^t_i,S^2(u^m_x))\br (S(u^l_y),u^z_s)
\br (S^2(u^x_k),u^j_z)\theta _{mn}F_1{}^s_aF_2{}^b_r\\
&\quad=\br (S^2(u^b_r),S(u^y_n))\br (u^t_i,S^2(u^m_x))\br (S(u^l_y),u^z_s)
\br (S^2(u^x_k),u^j_z)\theta _{mn}F_1{}^s_aF_2{}^r_t\\
&\quad=\br (S(u^b_r),u^y_n)\br (u^t_i,S^2(u^m_x))\br (S(u^l_y),S^2(u^s_a))
\br (S^2(u^x_k),u^j_z)\theta _{mn}F_1{}^z_sF_2{}^r_t\\
&\quad=\br (S(u^b_r),u^y_n)\br (S^2(u^r_t),S^2(u^m_x))\br (u^l_y,S(u^s_a))
\br (S^2(u^x_k),S^2(u^z_s))\theta _{mn}F_1{}^j_zF_2{}^t_i\\
&\quad=\br (u^r_t,u^m_x)\br (u^x_k,u^z_s)\br (S(u^b_r),u^y_n)
\br (u^l_y,S(u^s_a))\theta _{mn}F_1{}^j_zF_2{}^t_i\\
&\quad=g_{12}\bigl(\omega _{ij}\otA \sigma ^{-1}(\theta _{kl}\otA
\theta_{ab})\bigr).
\end{align*}
\end{bew}

Let us introduce the functional $f:\A \to \comp $ (see \cite{b-KS},
Proposition 10.3) defined by
$f(a)=\br (a_{(1)},S(a_{(2)}))$ and let $\bar{f}$ denote
the convolution inverse of $f$,
i.\,e.\ $\bar{f}(a)=\br (S^2(a_{(1)}),a_{(2)})$.

\begin{satz}\label{s-3bedmet}
Let $g$ be as in Proposition \ref{s-4bedmet}.
Then the bilinear map $g$ is $\sigma $-symmetric if and only if there are
complex numbers $c$ and $z$ such that $f(S(u^i_j))=z\bar{f}(u^i_j)$ and
\begin{gather*}
F_1{}^i_j=cG_2{}^i_kf(u^k_j)\qquad\text{and}\qquad
F_2{}^i_j=c^{-1}\bar{f}(u^i_k)G_1{}^k_j
\end{gather*}
for $i,j=1,\ldots,d$.
\end{satz}

\begin{bew}
Firstly let us suppose that $g\sigma =g$.
{}From the equation $g(\sigma (\omega _{ij}\otA \theta _{kl}))=
g(\omega _{ij}\otA \theta _{kl})$ we conclude that there is a nonzero complex
number $c'$ such that
\begin{gather*}\tag{$*$}
F_1{}^i_j=c'G_2{}^i_k\bar{f}(S(u^k_j))\quad\text{and}\quad
F_2{}^i_j=c'{}^{-1}f(S(u^i_k))G_1{}^k_j.
\end{gather*}
Further, $g(\sigma (\theta _{ij}\otA \omega _{kl}))=
g(\theta _{ij}\otA \omega _{kl})$ gives $G_1{}^i_j=cf(u^i_k)F_2{}^k_j$
and $G_2{}^i_j=c^{-1}F_1{}^i_k\bar{f}(u^k_j)$ for some nonzero
complex number $c$. Inserting this into ($*$) we obtain
$F_2{}^i_j=c'{}^{-1}cf(S(u^i_k))f(u^k_l)F_2{}^l_j$
for any $i,j$. Multiplying by $(F_2^{-1})^j_m\bar{f}(u^m_n)$ and summing up
over $j$ we obtain
$\bar{f}(u^i_n)=c'{}^{-1}cf(S(u^i_n))$. Inverting this equation we also get
$f(u^i_j)=c'c^{-1}\bar{f}(S(u^i_j))$.
Let us set $z=c'c^{-1}$. Then ($*$) gives the assertion.
The converse direction is an easy computation.
\end{bew}

\section{Contractions}
\label{sec-kontraktionen}

Let $\Gp $ and $\Gm $ be two bicovariant $\A $-bimodules over the
Hopf algebra $\A $.
Let $\Gaw{\tau }=\bigoplus_{k=0}^{\infty }\Gaw[k]{\tau }$,
$\tau \in \{+,-\}$ denote the external algebra for $\Ga{\tau }$ as
constructed by Woronowicz \cite{a-Woro2}. This means that there is an
antisymmetrizer $A_k: \Gat[k]{\tau }\to \Gat[k]{\tau }$ for each $k\geq 0$
($A_0=A_1=\id $) which is a
homomorphism of bicovariant bimodules and 
$\Gaw[k]{\tau }=\Gat[k]{\tau }/\ker A_k$.
Let us recall some properties of $A_k$.
Because of the general theory there are bimodule homomorphisms
$A_{i,j},B_{i,j}:\Gat[i+j]{\tau }\to \Gat[i+j]{\tau }$, $i,j\geq 0$
such that
\begin{gather}\label{eq-AB}
A_{i+j}=A_{i,j}(A_i\otA A_j),\quad
A_{i+j}=(A_i\otA A_j)B_{i,j}.
\end{gather}
In particular we have
\begin{align}\label{eq-aia}
A_i&=\prod_{k=0}^{i-1}(A_{i-k-1,1}\otA \id^{\otimes k})
=\prod_{k=0}^{i-1}(\id^{\otimes k}\otA A_{1,i-k-1})\\
\label{eq-aib}
A_i&=\prod_{k=0}^{i-1}(B_{k,1}\otA \id^{\otimes i-k-1})
=\prod_{k=0}^{i-1}(\id^{\otimes i-k-1}\otA B_{1,k}),
\end{align}
where
$A_{0,0}=A_{1,0}=A_{0,1}=\id $,
$B_{0,0}=B_{1,0}=B_{0,1}=\id $ and
\begin{align}\label{eq-a1i}
A_{1,i}&=\id -\sigma _{12}+\sigma _{23}\sigma _{12}-\ldots
+(-1)^i\sigma _{i,i+1}\cdots \sigma _{12},\\
\label{eq-ai1}
A_{i,1}&=\id -\sigma _{i,i+1}+\sigma _{i-1,i}\sigma _{i,i+1}-\ldots
+(-1)^i\sigma _{12}\cdots \sigma _{i,i+1},\\
\label{eq-b1i}
B_{1,i}&=\id -\sigma _{12}+\sigma _{12}\sigma _{23}-\ldots
+(-1)^i\sigma _{12}\cdots \sigma _{i,i+1},\\
\label{eq-bi1}
B_{i,1}&=\id -\sigma _{i,i+1}+\sigma _{i,i+1}\sigma _{i-1,i}-\ldots
+(-1)^i\sigma _{i,i+1}\cdots \sigma _{12}.
\end{align}
It is easy to see that
\begin{xalignat}{2}
\label{eq-a1rek}
A_{1,i}&=\id -(\id \otA A_{1,i-1})\sigma _{12},&
A_{i,1}&=\id -(A_{i-1,1}\otA \id )\sigma _{i,i+1},\\
\label{eq-b1rek}
B_{1,i}&=\id -\sigma _{12}(\id \otA B_{1,i-1}),&
B_{i,1}&=\id -\sigma _{i,i+1}(B_{i-1,1}\otA \id )
\end{xalignat}
for $i>0$.
One could also take (\ref{eq-a1i}) and (\ref{eq-aia}) for the definition
of $A_k$.

The preceding properties hold for any $\A $-bimodule isomorphism $\sigma $
which satisfies the braid relation. Therefore, replacing everywhere
$\sigma $ by $\sigma ^{-1}$ the above works as well.
In what follows we will use both kinds of operators and write
$A_k^{\tau }$, $A_{i,j}^{\tau }$ and $B_{i,j}^{\tau }$ whenever we are
dealing with $\sigma ^{\tau }$ ($\tau \in \{+,-\}$).

Let us introduce some operators in
$\mathrm{End}(\bigotimes _{i=1}^m\Ga{\tau _i})$,
$m\geq 1$ and $1\leq j,k\leq m$
(they can be associated to the permutations $(j,j+1,\ldots ,k)$,
$(k,k-1,\ldots ,j)$, $(1.m)(2,m-1)(3,m-2)\cdots $ and
$(1,2,\ldots ,j+1)(2,3,\ldots ,j+2)\cdots (k,k+1,\ldots j+k)$):
\begin{align}
\label{eq-Rto}
\sigma _\Rto{j}{k}^\pm &:=\sigma _{j,j+1}^\pm \sigma _{j+1,j+2}^\pm \cdots
\sigma _{k-1,k}^\pm \quad \text{for $j<k$},\quad
\sigma _\Rto{j}{k}^\pm =\id \quad \text{for $j\geq k$},\\
\label{eq-Lto}
\sigma _\Lto{j}{k}^\pm &:=\sigma _{k-1,k}^\pm \sigma _{k-2,k-1}^\pm \cdots
\sigma _{j,j+1}^\pm \quad \text{for $j<k$},\quad
\sigma _\Lto{j}{k}^\pm =\id \quad \text{for $j\geq k$},\\
\label{eq-lwortR}
\sigma _{(m)}^\pm &:=\sigma _\Lto{1}{1}^\pm \sigma _\Lto{1}{2}^\pm
\cdots \sigma _\Lto{1}{m}^\pm \quad \text{for $m\geq 1$},\quad
\sigma _{(0)}^\pm =\id,\\
\label{eq-bigsigR}
\sigma _{(j,k)}^\pm &:=\sigma _\Rto{k}{j+k}^\pm \sigma _\Rto{k-1}{j+k-1}^\pm
\cdots \sigma _\Rto{1}{j+1}^\pm \quad \text{for $m=j+k$}.
\end{align}
The verification of the following equations needs only braid group
techniques and is left to the reader. We have
\begin{gather}
\label{eq-lwortrek}
\sigma ^\pm _{(k)}=\sigma _\Rto{1}{k}(\sigma ^\pm _{(k-1)}\otA \id )=
\sigma ^\pm _\Lto{1}{k}(\id \otA \sigma ^\pm _{(k-1)}),\\
\label{eq-RLtoLR}
\sigma ^\pm _\Rto{1}{k}\sigma ^\pm _\Lto{1}{k-1}=
\sigma ^\pm _\Lto{1}{k}\sigma ^\pm _\Rto{2}{k}\\
\intertext{for $k\geq 2$ and}
\label{eq-Rtovert}
\sigma ^\pm _\Rto{1}{k}(b_{k-1}\otA \id)
=(\id \otA b_{k-1})\sigma ^\pm _\Rto{1}{k},\\
\label{eq-Ltovert}
\sigma ^\pm _\Lto{1}{k}(\id \otA b_{k-1})
=(b_{k-1} \otA \id )\sigma ^\pm _\Lto{1}{k},\\
\label{eq-bigsigL}
\sigma ^\pm _{(j,k)}=\sigma ^\pm _\Lto{1}{k+1}\sigma ^\pm _\Lto{2}{k+2}\cdots
\sigma ^\pm _\Lto{j}{k+j}
\end{gather}
for $k\geq 1$, where $b_k$ is an arbitrary expression of the complex
algebra generated by $\sigma _{12},\ldots ,\sigma _{k-1,k}$ and their
inverses.
Observe that $A^\pm _k\sigma ^\mp _{(k)}=\sigma ^\mp _{(k)}A^\pm _k=
(-1)^{k(k-1)/2}A^\mp _k$ (see \cite{a-Woro2}, p.\ 157).
Hence, in particular have $\ker A^+ _k=\ker A^- _k$.

Now let $g$ be a $\sigma $-metric of the pair $(\Gp ,\Gm )$.
The next formulas follow from the fourth condition on the $\sigma $-metric
by induction over k:
\begin{gather}\label{eq-gRRgLL}
g_{12}\sigma ^\pm _\Rto{2}{k} \sigma ^\pm _\Rto{1}{k-1}=g_{k-1,k},\quad
g_{k-1,k}\sigma ^\pm _\Lto{1}{k-1} \sigma ^\pm _\Lto{2}{k}=g_{12}\quad
\text{for $k\geq 2$.}
\end{gather}

Next we define \textit{contractions} $\ctrpm{\cdot }{\cdot }:
\Gat[k]{\tau }\tet \Gat[l]{-\tau }\to \Gat[|k-l|]{\tau '}$,
$\tau \in \{+,-\}$, $\tau '=\tau$ for $k\geq l$, otherwise $\tau '=-\tau $,
by
\begin{equation}\label{eq-ctr}
\begin{aligned}
\ctrpm{\xi }{\xi '}&:=\tilde{g}(B^\pm _{k-l,l}\xi \mt A^\pm _l\xi ')\quad
\text{for $k\geq l$,}\\
\ctrpm{\xi }{\xi '}&:=\tilde{g}(A^\pm _k\xi \mt B^\pm _{k,l-k}\xi ')\quad
\text{for $k<l$.}
\end{aligned}
\end{equation}
This maps are homomorphisms of $\A $-bimodules and inherit all covariance
properties of $g$.
If both $k$ and $l$ are less than two, then the contraction doesn't depend
on the sign $\pm $ and we sometimes omit it:
$\ctrp{\xi }{\xi '}=\ctrm{\xi }{\xi '}=:\ctr{\xi }{\xi '}$.

Next we prove a generalisation of Lemma \ref{l-tgsig}.

\begin{lemma}\label{l-tgant}
Let $g$ be a $\sigma $-metric of the pair $(\Gp ,\Gm )$ and let
$\tilde{g}$ be the
map defined by (\ref{eq-tildeg}). Then we have for all nonnegative
integers $i,j,k,l$, $1\leq i+j\leq k,l$,
\begin{align*}
\tilde{g}\circ \bigl( (\id ^{\otimes k-i-j}\otA A^\pm _i \otA \id ^{\otimes j})
\mt \id ^{\otimes l}\bigr) &=
\tilde{g}\circ \bigl( \id ^{\otimes k}\mt (\id ^{\otimes j}\otA A^\pm _i
\otA \id ^{\otimes l-i-j})\bigr) .
\end{align*}
\end{lemma}

\begin{bew}
Using Lemma \ref{l-tgsig} one checks that
\begin{align}
\tilde{g}\circ (\sigma ^\pm _\RLto{k+1-t'}{k+1-t}\mt \id ^{\otimes l})&=
\tilde{g}\circ (\id ^{\otimes k}\mt \sigma ^\pm _\RLto{t}{t'})
\end{align}
for $1\leq t\leq t'\leq k,l$.
{}From this and equations (\ref{eq-bi1}) and (\ref{eq-a1i}) we obtain
\begin{align*}
&\tilde{g}\circ \bigl( (\id ^{\otimes k-r-s}\otA A^\pm _{1,s-1}\otA
\id ^{\otimes r}) \mt \id ^{\otimes l}\bigr) =\\
&\quad =
\tilde{g}\circ \left( \id ^{\otimes k-r-s}\otA \left(\sum _{t=1}^s
(-1)^{t+1}\sigma ^\pm _\Lto{1}{t}\right) \otA \id ^{\otimes r}\mt
\id ^{\otimes l}\right) \\
&\quad =\tilde{g}\circ \left( \left( \sum _{t=1}^s (-1)^{t+1}
\sigma ^\pm _\Lto{k+1-r-s}{k-r-s+t}\right) \mt \id ^{\otimes l}\right) \\
&\quad =\tilde{g}\circ \left( \id ^{\otimes k}\mt \left(
\sum _{t=1}^s(-1)^{t+1}\sigma ^\pm _\Lto{r+s-t+1}{r+s}\right) \right) \\
&\quad =\tilde{g}\circ \left( \id ^{\otimes k}\mt \id ^{\otimes r}\otA
\left( \sum _{t=1}^s(-1)^{t+1}\sigma ^\pm _\Lto{s-t+1}{s}\right) \otA
\id ^{\otimes l-r-s}\right) \\
&\quad =\tilde{g}\circ \bigl( \id ^{\otimes k}\mt (\id ^{\otimes r}\otA
B^\pm _{s-1,1}\otA \id ^{\otimes l-r-s})\bigr)
\end{align*}
for all $r,s$ with $0\leq r$, $1\leq s$, $r+s\leq k$, $r+s\leq l$.
Using this result together with (\ref{eq-aia}) and (\ref{eq-aib})
similar computations give the assertion of the lemma.
\end{bew}

\begin{lemma}\label{l-tgbigsig}
Let $\tilde{g}$ be as in Lemma \ref{l-tgant} and $\xi _k\in \Gat[k]{\tau }$,
$\xi '_l\in \Gat[l]{-\tau }$, $\tau \in \{+,-\}$, $k,l\geq 0$. Then
$\sigma ^\pm _{(k,l)}(\sigma ^\pm _{(k)}(\xi _k)\mt
\sigma ^\pm _{(l)}(\xi '_l))$ is an element of
$\Gat[l]{-\tau }\tet \Gat[k]{\tau }$ and the equation
\begin{align}\label{eq-tgbigsig}
\tilde{g}\bigl( \sigma ^\pm _{(k,l)}(\sigma ^\pm _{(k)}(\xi _k)\mt
\sigma ^\pm _{(l)}(\xi '_l)) \bigr) &=
\sigma ^\pm _{(|k-l|)}(\tilde{g}(\xi _k\mt \xi '_l))
\end{align}
holds.
\end{lemma}

\begin{bew}
We prove the case $k\geq l$ by induction on $l$. Then the assertion follows
also for $k<l$ because of the formulas
$\sigma ^\pm _{(i)}\sigma ^\mp _{(i)}=\id$,
$\sigma ^\pm _{(i,j)}(A^\pm _{i}\mt \id ^{\otimes j})=
(\id ^{\otimes j}\mt A^\pm _{i})\sigma ^\pm _{(i,j)}$
(see also (\ref{eq-bigsigR}) and (\ref{eq-lwortrek})) and
$\sigma ^\pm _{(i,j)}\sigma ^\mp _{(j,i)}=\id$ for all $i,j\geq 0$.

If $l=0$ then $\sigma ^\pm _{(l)}=\sigma ^\pm _{(k,l)}=\id$, hence
the left hand side of (\ref{eq-tgbigsig}) is equal to
$\tilde{g}(\sigma ^\pm _{(k)}(\xi _k)\mt \xi '_l)=
\sigma ^\pm _{(k)}(\xi _k)\xi '_l$.
For the right hand side we obtain
$\sigma ^\pm _{(k-l)}\tilde{g}(\xi _k\mt \xi '_l)
=\sigma ^\pm _{(k-l)}(\xi _k\xi '_l)$. Since $\xi '_l\in \A$ and
$\sigma ^\pm _{(k-l)}$ is a homomorphism of $\A $-bimodules, the assertion
of the lemma is valid.

Suppose that (\ref{eq-tgbigsig}) holds for an $l\in\mathbb{N}_0$, $l\leq k$.
Consider the map
\[\tilde{g}\sigma ^\pm _{(k+1,l+1)}
(\sigma ^\pm _{(k+1)}\mt \sigma ^\pm _{(l+1)}):
\Gat[k+1]{\tau }\tet \Gat[l+1]{-\tau }\to \Gat[k-l]{\tau }.
\]
We compute
\begin{align*}
&\tilde{g}\sigma ^\pm _{(k+1,l+1)}
(\sigma ^\pm _{(k+1)}\otA \sigma ^\pm _{(l+1)})\\
&\quad
=\tilde{g}g_{l+1,l+2}\sigma ^\pm _\Rto{l+1}{k+l+2}
(\sigma ^\pm _{(k+1,l)}\otA \id)(\sigma ^\pm _\Rto{1}{k+1}(\sigma ^\pm _{(k)}
\otA \id )\otA \sigma ^\pm _{(l+1)})\\
&\quad
=\tilde{g}g_{l+1,l+2}\sigma ^\pm _\Rto{l+2}{k+l+2}
\sigma ^\pm _\Rto{k+1}{k+l+1}
(\sigma ^\pm _{(k+1,l)}\otA \id)(\sigma ^\pm _{(k)} \otA \id \otA
\sigma ^\pm _{(l+1)})\\
&\quad
=\tilde{g}g_{k+l+1,k+l+2}
(\sigma ^\pm _{(k+1,l)}\otA \id)(\sigma ^\pm _{(k)} \otA \id \otA
\sigma ^\pm _{(l+1)})\\
&\quad
=\tilde{g}g_{k+l+1,k+l+2}(\sigma ^\pm _{(k,l)}\otA \id ^{\otimes 2})
\sigma ^\pm _\Lto{k+1}{k+l+1}\sigma ^\pm _\Lto{k+2}{k+l+2}
(\sigma ^\pm _{(k)} \otA \id ^{\otimes 2}\otA \sigma ^\pm _{(l)})\\
&\quad
=\tilde{g}\sigma ^\pm _{(k,l)}g_{k+l+1,k+l+2}
\sigma ^\pm _\Lto{k+1}{k+l+1}\sigma ^\pm _\Lto{k+2}{k+l+2}
(\sigma ^\pm _{(k)} \otA \id ^{\otimes 2}\otA \sigma ^\pm _{(l)})\\
&\quad
=\tilde{g}\sigma ^\pm _{(k,l)}g_{k+1,k+2}
(\sigma ^\pm _{(k)} \otA \id ^{\otimes 2}\otA \sigma ^\pm _{(l)})\\
&\quad
=\tilde{g}\sigma ^\pm _{(k,l)}
(\sigma ^\pm _{(k)} \otA g_{12}\otA \sigma ^\pm _{(l)})
=\tilde{g}g_{k+1,k+2}
=\tilde{g}
\end{align*}
where we used the following formulas:
(\ref{eq-bigsigR}) and (\ref{eq-lwortrek}) in the first equation,
the $\sigma $-symmetry of the $\sigma $-metric, (\ref{eq-bigsigR}) and
(\ref{eq-Rtovert}) in the second,
(\ref{eq-gRRgLL}) in the third,
(\ref{eq-bigsigL}) and (\ref{eq-lwortrek}) in the fourth,
(\ref{eq-gRRgLL}) in the sixth,
the induction assumption in the eighth and the recursive definition
of $\tilde{g}$ in the last equation.
\end{bew}

An important consequence of Lemma \ref{l-tgant} is the
possibility to extend the definition of our contractions
$\ctrpm{\cdot }{\cdot }$ to a map
$\ctrpm{\cdot }{\cdot }:
\Gaw[k]{\tau }\tet \Gaw[l]{-\tau }\to \Gaw[|k-l|]{\tau '}$,
$\tau \in \{+,-\}$, $\tau '=\tau$ for $k\geq l$, otherwise $\tau '=-\tau $.
To see this, we treat the case $k \geq l$.
Let $\xi '_k\in \Gat[k]{\tau }$ and $\xi _l\in \Gat[l]{-\tau }$,
$\tau \in \{+,-\}$. Firstly, let
$\xi _l$ be a symmetric $l$-form, i.\,e.\ $A^\pm _l(\xi _l)=0$.
Then, by definition,
\begin{gather*}
\ctrpm{\xi '_k}{\xi _l}=
\tilde{g}(B^\pm _{k-l,l}\xi '_k\mt A^\pm _l\xi _l)=0.
\end{gather*}
On the other hand, if
$\xi '_k$ is a symmetric $k$-form, i.\,e.\ $A^\pm _k(\xi '_k)=0$,
then we conclude
\begin{gather*}
A^\pm _{k-l}\ctrpm{\xi '_k}{\xi _l}=
A^\pm _{k-l}\tilde{g}(B^\pm _{k-l,l}\xi '_k\mt A^\pm _l\xi _l).
\end{gather*}
Applying Lemma \ref{l-tgant} this is equal to
\begin{gather*}
A^\pm _{k-l}\tilde{g}((\id ^{\otimes k-l}\otA A^\pm _l)B^\pm _{k-l,l}\xi '_k
\mt \xi _l)=
\tilde{g}((A^\pm _{k-l}\otA A^\pm _l)B^\pm _{k-l,l}\xi '_k\mt \xi _l).
\end{gather*}
Now formula (\ref{eq-AB}) insures that the latter expression is zero.
Hence $\ctrpm{\xi '_k}{\xi _l}$ is symmetric.
In the case $k<l$ similar reasoning gives the desired result.

\begin{bem}
In view of Lemma \ref{l-symmet} and Proposition \ref{s-metna} we should
also consider the contractions for $k=l$ (composed with the Haar functional,
see in Section \ref{sec-Laplace}) as a kind of higher rank $\sigma $-metric.
\end{bem}

\begin{lemma}\label{l-ctrrek}
For $\xi _i\in \Gaw[k_i]{\tau _i}$, $i=0,1,2$, $\tau _1=\tau _2=-\tau _0$,
$k_1+k_2\leq k_0$ the contractions satisfy the following relations:\\
(i) $\ctrpm{\xi _1}{\ctrpm{\xi _2}{\xi _0}}=
\ctrpm{\xi _1 \wedge \xi _2}{\xi _0}$ and
$\ctrpm{\ctrpm{\xi _0}{\xi _1}}{\xi _2}=
\ctrpm{\xi _0}{\xi _1\wedge \xi _2}$,\\
(ii) $\ctrpm{\xi _1}{\ctrpm{\xi _0}{\xi _2}}=
\ctrpm{\ctrpm{\xi _1}{\xi _0}}{\xi _2}$.
\end{lemma}

\begin{bew}
{}From Lemma \ref{l-tgant} and formula (\ref{eq-AB}) we conclude that
$A^\pm _{k-l}\ctrpm{\zeta '_l}{\zeta ''_k}=\tilde{g}(\zeta '_l\mt
A^\pm _k\zeta ''_k)$ for $k\geq l$, $\zeta ''_k\in \Gaw[k]{\tau }$,
$\zeta '_l\in \Gaw[l]{-\tau }$, $\tau \in \{+,-\}$.
Then for the first equation of (i) and representants
$\zeta _i\in \Gat[k_i]{\tau _i}$ of $\xi _i$, $i=0,1,2$ we compute
\begin{align*}
A^\pm _{k_0-k_1-k_2}(\ctrpm{\zeta _1}{\ctrpm{\zeta _2}{\zeta _0}})&=
\tilde{g}(\zeta _1 \mt A^\pm _{k_0-k_2}(\ctrpm{\zeta _2}{\zeta _0}))\\
&=\tilde{g}(\zeta _1 \mt \tilde{g}(\zeta _2 \mt A^\pm _{k_0}\zeta _0))
=\tilde{g}(\zeta _1 \otA \zeta _2 \mt A^\pm _{k_0}\zeta _0)\\
&=A^\pm _{k_0-k_1-k_2}(\ctrpm{\zeta _1 \otA \zeta _2}{\zeta _0}).
\end{align*}
The second equation can be proved similarly.

To prove (ii) we use the same arguments. For the left hand side we obtain
\begin{align*}
A^\pm _{k_0-k_1-k_2}\ctrpm{\zeta _1}{\ctrpm{\zeta _0}{\zeta _2}}&=
\tilde{g}(\zeta _1 \mt A^\pm _{k_0-k_2}\ctrpm{\zeta _0}{\zeta _2})\\
&=\tilde{g}(\zeta _1 \mt \tilde{g}(A^\pm _{k_0}\zeta _0 \mt \zeta _2))
\end{align*}
and for the right hand side
\begin{align*}
A^\pm _{k_0-k_1-k_2}\ctrpm{\ctrpm{\zeta _1}{\zeta _0}}{\zeta _2}&=
\tilde{g}(A^\pm _{k_0-k_1}\ctrpm{\zeta _1}{\zeta _0}\mt \zeta _2)\\
&=\tilde{g}(\tilde{g}(\zeta _1 \mt A^\pm _{k_0}\zeta _0) \mt \zeta _2).
\end{align*}
But both last expressions are equal because of the definition
of $\tilde{g}$ and since $k_1+k_2\leq k_0$.
\end{bew}

The following lemma contains some recursion formulas which are useful
in order to compute contractions.

\begin{lemma}
For any $\xi _k\in \Gaw[k]{\tau }$, $\xi '_k\in \Gaw[k]{-\tau }$,
$\rho _1\in \Ga{\tau}$, $\rho _2\in \Ga{-\tau }$, $k\geq 1$,
$\tau \in \{+,-\}$ the equations
\begin{align}
\label{eq-ctrrekr}
\ctrpm{\xi _k\wedge \rho _1}{\rho _2}&=\xi _k\ctrpm{\rho _1}{\rho _2}-
\ctrpm{\xi _k}{\rho ^\mp _{(1)}}\wedge \rho ^\mp _{(2)}\quad \text{and}\\
\label{eq-ctrrekl}
\ctrpm{\rho _1}{\rho _2\wedge \xi '_k}&=\ctrpm{\rho _1}{\rho _2}\xi '_k
-\rho ^\mp _{(1)}\wedge \ctrpm{\rho ^\mp _{(2)}}{\xi '_k}
\end{align}
hold, where
$\sigma ^\mp (\rho _1\otA \rho _2)=\rho ^\mp _{(1)}\otA \rho ^\mp _{(2)}
\in \Ga{-\tau }\otA\Ga{\tau }$.
\end{lemma}

\begin{bew}
For $k=1$ the left hand side of the first equation reads as
\begin{align*}
\ctrpm{\xi _1\wedge \rho _1}{\rho _2}&=
\tilde{g}\bigl( (\xi _1\otA \rho _1-\sigma ^\pm (\xi _1\otA \rho _1))\mt
\rho _2\bigr)\\
&=g_{23}(\xi _1\otA \rho _1\otA \rho _2)-
g_{12}\sigma ^\mp _{23}(\xi _1\otA \rho _1\otA \rho _2)
\end{align*}
because of the fourth condition of Definition \ref{d-metrik} on the
$\sigma $-metric $g$.
Further, if $k>1$ we then use (\ref{eq-b1rek}) to conclude in a similar manner
that
\begin{align*}
&\ctrpm{\xi _k\wedge \rho _1}{\rho _2}=
\tilde{g}(B_{k,1}(\xi _k\otA \rho _1)\mt \rho _2)\\
&\quad =g_{k+1,k+2}(\xi _k\otA \rho _1\otA \rho _2
-\sigma ^\pm _{k,k+1}(B_{k-1,1}\otA \id ^{\otimes 2})(\xi _k\otA
\rho _1 \otA \rho _2))\\
&\quad =\xi _k\otA g(\rho _1\otA \rho _2)
-g_{k,k+1}\sigma ^\mp _{k+1,k+2}(B_{k-1,1}\otA \id ^{\otimes 2})
(\xi _k\otA \rho _1 \otA \rho _2)\\
&\quad =\xi _k\otA g(\rho _1\otA \rho _2)
-g_{k,k+1}(B_{k-1,1}\otA \sigma ^\mp )
(\xi _k\otA \rho _1 \otA \rho _2)\\
&\quad =\xi _k\ctrpm{\rho _1}{\rho _2}-
\ctrpm{\xi _k}{\rho ^\mp _{(1)}}\wedge \rho ^\mp _{(2)}.
\end{align*}
The proof of the second equation of the lemma is analogous.
\end{bew}

\begin{lemma}\label{l-symmet}
For fixed $\tau \in \{+,-\},k\geq 1$ let $\rho _k\in \Gaw[k]{\tau },
\rho '_k\in\Gaw[k]{-\tau }$ and
$\sigma ^\pm _{(k,k)}(\rho _k\otA \rho '_k)
=\rho '{}^k_{(1)}\otA \rho ^k_{(2)}$.
Then the contractions $\ctrpm{\cdot }{\cdot }$ satisfy the equations
\begin{align}
\ctrpm{\rho '{}^k_{(1)}}{\rho ^k_{(2)}}&=
\ctrpm{\rho _k}{(\sigma ^\mp _{(k)})^2(\rho '_k)}=
\ctrpm{(\sigma ^\mp _{(k)})^2(\rho _k)}{\rho '_k}.
\end{align}
\end{lemma}

\begin{bew}
The definition (\ref{eq-ctr}) of $\ctrpm{\cdot }{\cdot }$ gives
$\ctrpm{\rho '{}^k_{(1)}}{\rho ^k_{(2)}}
=\tilde{g}(\rho '{}^k_{(1)}\mt A^\pm _k\rho ^k_{(2)})$.
Since $(\id ^{\otimes k}\otA A^\pm _k)\sigma ^\pm _{(k,k)}
=\sigma ^\pm _{(k,k)}(A^\pm _k\otA \id ^{\otimes k})$
(see (\ref{eq-bigsigL}) and (\ref{eq-Ltovert}))
and $A^\pm _k=(-1)^{k(k-1)/2}A^\mp _k\sigma ^\pm _{(k)}$,
we conclude
\begin{align*}
\ctrpm{\rho '{}^k_{(1)}}{\rho ^k_{(2)}}&=
\tilde{g}(\id ^{\otimes k}\otA A^\pm _k)\sigma ^\pm _{(k,k)}
(\rho _k\otA \rho '_k)\\
&=\tilde{g}\sigma ^\pm _{(k,k)}(A^\pm _k\otA \id ^{\otimes k})
(\rho _k\otA \rho '_k)\\
&=(-1)^{k(k-1)/2}\tilde{g}\sigma ^\pm _{(k,k)}(A^\mp _k\sigma ^\pm _{(k)}
\rho _k\otA \rho '_k)\\
&=(-1)^{k(k-1)/2}\tilde{g}\sigma ^\pm _{(k,k)}(\sigma ^\pm _{(k)}
\rho _k\otA A^\mp _{(k)}\rho '_k)
\end{align*}
by Lemma \ref{l-tgant}. Inserting
$(-1)^{k(k-1)/2}A^\mp _k=\sigma ^\pm _{(k)}A^\pm _k(\sigma ^\mp _{(k)})^2$
and applying Lemma \ref{l-tgbigsig} we obtain
\begin{align*}
\ctrpm{\rho '{}^k_{(1)}}{\rho ^k_{(2)}}&=\tilde{g}\sigma ^\pm _{(k,k)}
\bigl( \sigma ^\pm _{(k)}\rho _k\mt \sigma ^\pm _{(k)}A^\pm _k
(\sigma ^\mp _{(k)})^2\rho '_k\bigr)\\
&=\tilde{g}(\rho _k\mt A^\pm _k(\sigma ^\mp _{(k)})^2\rho '_k)
=\ctrpm{\rho _k}{(\sigma ^\mp _{(k)})^2\rho '_k}.
\end{align*}
The second equation follows similarly.
\end{bew}

Finally, we should say something about the nondegeneracy of
$\ctrpm{\cdot }{\cdot }$ as a $\sigma $-metric.

\begin{satz}\label{s-metna}
The maps $\ctrpm{\cdot }{\cdot }:\Gaw[k]{\tau }\tet \Gaw[k]{-\tau }\to \A $,
$\tau \in \{+,-\}$, $k\geq 1$ and their restrictions to
$\linv{\Gaw[k]{\tau }}\mt \linv{\Gaw[k]{-\tau }}$ are nondegenerate.
\end{satz}

\begin{bew}
Firstly we show that $\tilde{g}:\Gat[k]{\tau }\tet \Gat[k]{-\tau }\to \A $
and its restriction to $\linv{\Gat[k]{\tau }}\mt \linv{\Gat[k]{-\tau }}$
are nondegenerate.

For $k=1$ this assertion is true, since $g$ is
nondegenerate by Definition \ref{d-metrik} and $\tilde{g}=g$.
Suppose that it is valid for some $k\geq 1$ and let
$\xi _{k+1}\in \Gat[k+1]{\tau }$. Then there are finitely many $k$-forms
$\xi ^i\in\Gat[k]{\tau }$ and linearly independent
1-forms $\rho _i\in \linv{\Ga{\tau }}$ such that
$\xi _{k+1}=\sum _i\rho _i\otA \xi ^i$.
Suppose that $\tilde{g}(\xi _{k+1}\mt (\xi '_k\otA \rho '))=0$
for any $\xi '_k\in\linv{\Gat[k]{-\tau }}$ and $\rho '\in\linv{\Ga{-\tau }}$.
Hence by definition of $\tilde{g}$,
$\tilde{g}((\rho _i\otA \xi ^i)\mt (\xi '_k\otA \rho '))=
g(\rho _i\tilde{g}(\xi ^i\mt \xi '_k)\mt \rho ')=0$ for any
$\rho '\in \linv{\Ga{-\tau }}$. Since $g$ is a homomorphism of right
$\A $-modules, the latter is also true for any $\rho '\in\Ga{-\tau }$.
Applying the nondegeneracy of $g$ we conclude that
$\rho _i\tilde{g}(\xi ^i\mt \xi '_k)=0$
and since the 1-forms $\rho _i\in \linv{\Ga{\tau }}$ are linearly
independent we obtain $\tilde{g}(\xi ^i\mt \xi '_k)=0$ for any
$\xi '_k\in \linv{\Gat[k]{-\tau }}$.
Now we use that $\tilde{g}$ is a homomorphism of right $\A $-modules and get
$\tilde{g}(\xi ^i\mt \xi '_k)=0$ for any $\xi '_k\in \Gat[k]{-\tau }$.
Then the induction assumption gives $\xi ^i=0$ and hence
$\xi _{k+1}=\rho _i\otA \xi ^i=0$.

Now we prove the assertion of the proposition. Let
$\xi \in \Gaw[k]{\tau }$, $\xi _0\in \Gat[k]{\tau }$ be
a representant of $\xi $, and let us assume that $\ctrpm{\xi }{\xi '_k}=0$
for any $\xi '_k\in \linv{\Gaw[k]{-\tau }}$. This means
$\tilde{g}(A_k^\pm \xi _0\mt \xi '_k)=0$ for any
$\xi '_k\in \linv{\Gaw[k]{-\tau }}$. Since $\tilde{g}$ is a homomorphism
of right $\A $-modules, the latter is true for any
$\xi '_k\in \Gaw[k]{-\tau }$. In the first part of the proof we have shown
that $A_k^\pm \xi _0=0$. Hence $\xi _0$ is a symmetric $k$-form,
so that $\xi =0$.

Nondegeneracy in the second component of $\ctrpm{\cdot }{\cdot }$ can be
proved similarly.
\end{bew}

\begin{folg}\label{f-dimpmgl}
Let $g$ be a left-covariant $\sigma $-metric of the pair
$(\Gp ,\Gm )$. Then
for any $k\geq 0$ we have $\dim \linv{\Gaw[k]{+}}=\dim \linv{\Gaw[k]{-}}$.
\end{folg}

\section{Hodge Operators}
\label{sec-Hodge}

In this section we assume that
\begin{itemize}
\item[(I)]
the only one-dimensional corepresentation of the Hopf algebra
$\A $ is $1$ and
\item[(II)]
there exists a nonzero differential form
$\omega _0^{\tau }\in \linv{\Gaw[n]{\tau }}$ for some $n\in \mbb{Z}$
and $\tau \in \{+,-\}$ such that
$\omega _0^{\tau }\wedge \rho =0$ for all $\rho \in \Ga{\tau }$.
\end{itemize}
The latter is in particular fulfilled if one of the vector spaces
$\linv{\Gaw{+}}$, $\linv{\Gaw{-}}$ is finite dimensional.
Let us fix a triple $(n_0,\tau _0,\omega _0^{\tau _0})$ as in (II)
such that for any other triple $(n_1,\tau _1,\omega _1^{\tau _1} )$
having the same property we have $n_1\geq n_0$.

After proving some statements we will show that both $+$ and $-$
can occur as the value of $\tau _0$ and for a given left-covariant
$\sigma $-metric
$g$ of the pair $(\Gp ,\Gm )$, $\omega _0^\pm $ can be taken biinvariant
and in such a manner that
\begin{gather}\label{eq-ctrnorm}
\ctrpm{\omega _0^+}{\omega _0^-}=
\ctrpm{\omega _0^-}{\omega _0^+}=1.
\end{gather}
Then we also will assume this on $\omega _0^+$ and $\omega _0^-$.

Let $g$ be a (not necessarily left-covariant) $\sigma $-metric of the pair
$(\Gp ,\Gm )$.

\begin{satz}\label{s-hfumfr}
For any $\xi _k\in \Gaw[k]{-\tau _0}$, $\xi '_l\in \Gaw[l]{\tau _0}$,
$0\leq l\leq k\leq n_0$, we have
\begin{gather}\label{eq-hfumfr}
\ctrpm{\omega _0^{\tau _0}}{\xi _k}\wedge \xi '_l=
\ctrpm{\omega _0^{\tau _0}}{\ctrmp{\xi _k}{\xi '_l}}.
\end{gather}
\end{satz}

\begin{bew}
For $l=0$ the assertion follows from the right $\A $-linearity of
$\tilde{g}$.
Let us examine first the case $k=l=1$.
Inserting $\tau =-\tau _0$ and $\xi _k=\omega _0^{\tau _0}$ into
(\ref{eq-ctrrekr}) and using the condition on $\omega _0^{\tau _0}$
we obtain
$0=\omega _0^{\tau _0}\ctrpm{\rho _1}{\rho _2}-
\ctrpm{\omega _0^{\tau _0}}{\rho _{(1)}^\mp }\wedge \rho _{(2)}^\mp $
for any $\rho _1\in \Ga{\tau _0}$ and $\rho _2\in \Ga{-\tau _0}$, where
$\rho _{(1)}^\mp \otA \rho _{(2)}^\mp =\sigma ^\mp (\rho _1\otA \rho _2)$.
Now we insert $\sigma ^\pm (\xi _1\otA \xi '_1)$ for $\rho _1\otA \rho _2$
and obtain the desired result by the $\sigma $-symmetry of the
$\sigma $-metric $g$.

Secondly we prove the proposition for $1=l\leq k\leq n_0$ by induction on
$k$. The first step for this is already done. Suppose now that the assertion
is true for a $k<n_0$ and let $\xi _k\in \Gaw[k]{-\tau _0}$,
$\rho _1\in \Ga{\tau _0}$ and $\rho _2\in \Ga{-\tau _0}$.
By (\ref{eq-ctrrekr}) we obtain
\begin{align*}\tag{$*$}
\ctrpm{\ctrpm{\omega _0^{\tau _0}}{\xi _k}\wedge \rho _1}{\rho _2}&=
\ctrpm{\omega _0^{\tau _0}}{\xi _k}\ctrpm{\rho _1}{\rho _2}-
\ctrpm{\ctrpm{\omega _0^{\tau _0}}{\xi _k}}{\rho ^\mp _{(1)}}\wedge
\rho ^\mp _{(2)}
\end{align*}
where
$\rho _{(1)}^\mp \otA \rho _{(2)}^\mp =\sigma ^\mp (\rho _1\otA \rho _2)$.
The induction assumption and the second equation of Lemma \ref{l-ctrrek}(i)
assure that the left hand side of the latter equation is equal to
\begin{align*}
\ctrpm{\ctrpm{\omega _0^{\tau _0}}{\ctrmp{\xi _k}{\rho _1}}}{\rho _2}&=
\ctrpm{\omega _0^{\tau _0}}{\ctrmp{\xi _k}{\rho _1}\wedge \rho _2}.
\end{align*}
Moving this to the right hand side and the second term of the
right hand side of ($*$) to the left we get
\begin{align*}
\ctrpm{\ctrpm{\omega _0^{\tau _0}}{\xi _k}}{\rho ^\mp _{(1)}}\wedge
\rho ^\mp _{(2)}&=
\ctrpm{\omega _0^{\tau _0}}{\xi _k\ctrmp{\rho _1}{\rho _2}}-
\ctrpm{\omega _0^{\tau _0}}{\ctrmp{\xi _k}{\rho _1}\wedge \rho _2},
\end{align*}
where we used
the right $\A $-linearity of the contraction and the relation
$\ctrp{\rho _1}{\rho _2}=\ctrm{\rho _1}{\rho _2}$.
Now we take arbitrary elements $\xi '_1\in \Ga{\tau _0}$,
$\xi ''_1\in \Ga{-\tau _0}$.
We insert $\sigma ^\pm (\xi ''_1\otA \xi '_1)$ for $\rho _1\otA \rho _2$
in the above formula and use Lemma \ref{l-ctrrek}.(i) (on the left hand side),
the $\sigma $-symmetry of $g$ (in the first term of the right hand side)
and (\ref{eq-ctrrekr}) (on the right hand side of the latter equation).
In this manner we obtain
\begin{align*}
\ctrpm{\omega _0^{\tau _0}}{\xi _k \wedge \xi ''_1}\wedge \xi '_1&=
\ctrpm{\omega _0^{\tau _0}}{\ctrmp{\xi _k \wedge \xi ''_1}{\xi '_1}}.
\end{align*}
Hence the assertion of the proposition is true for $k+1$.

Suppose now that the assertion of the proposition is
valid for a fixed $l<n_0$ and for all $k>l$.
For $l=1$ this is true. Then for arbitrary $\xi ''\in \Ga{\tau _0}$
we apply (\ref{eq-hfumfr}) twice and conclude
\begin{align*}
\ctrpm{\omega _0^{\tau _0}}{\xi _k}\wedge \xi '_l\wedge \xi ''_l=
\ctrpm{\omega _0^{\tau _0}}{\ctrmp{\xi _k}{\xi '_l}}\wedge \xi ''_l=
\ctrpm{\omega _0^{\tau _0}}{\ctrmp{\ctrmp{\xi _k}{\xi '_l}}{\xi ''_l}}.
\end{align*}
Applying now Lemma \ref{l-ctrrek}.(i), we get (\ref{eq-hfumfr})
for $l+1$.
\end{bew}

{}From now on let $g$ be a left-covariant $\sigma $-metric of the pair
$(\Gp ,\Gm )$.
A very important consequence of Proposition \ref{s-hfumfr} is the following.

\begin{thm}\label{t-eindhf}
If there is a left-covariant $\sigma $-metric $g$ of the pair $(\Gp ,\Gm )$
then there exists a natural number $n_0$ such that
$\dim \linv{\Gaw[n_0]{\tau }}=1$ for $\tau \in \{+,-\}$ and all $k$-forms
$\xi _k\in \Gaw[k]{\tau }$, $k>n_0$ vanish.
\end{thm}

\begin{bew}
Since the $\sigma $-metric $\ctrpm{\cdot}{\cdot}$ is nondegenerate by
Proposition \ref{s-metna} and left-covariant
there is a left-invariant $n_0$-form $\xi _{n_0}\in \Gaw[n_0]{-\tau _0}$
such that $\ctrp{\omega _0^{\tau _0}}{\xi _{n_0}}=1$. Inserting
an arbitrary $\xi '_l$, $l=n_0$ into (\ref{eq-hfumfr})
we obtain $\xi '_{n_0}=\ctrp{\omega _0^{\tau _0}}{\xi _{n_0}}\xi '_{n_0}=
\omega _0^{\tau _0}\ctrm{\xi _{n_0}}{\xi '_{n_0}}$.
Hence we get $\Gaw[n_0]{\tau _0}=\omega _0^{\tau _0}\cdot \A $.
Since $\Gaw[k]{\tau _0}=\Gaw[n_0]{\tau _0}\wedge \Gaw[k-n_0]{\tau _0}
=\omega _0^{\tau _0}\wedge \Gaw[k-n_0]{\tau _0}$
for any $k>n_0$, we obtain $\Gaw[k]{\tau _0}={0}$.
The same assertion for $-\tau _0$ follows from Corollary \ref{f-dimpmgl}.
\end{bew}

\begin{bem}
In the proofs of Proposition \ref{s-hfumfr} and Theorem \ref{t-eindhf}
the assumption that there is only one one-dimensional corepresentation of
$\A $ was not used.
\end{bem}

\begin{folg}\label{f-eindhf}
Let $\A $ be an arbitrary Hopf algebra over the complex field with invertible
antipode.
Let $\Gp $ and $\Gm $ be bicovariant $\A $-bimodules and $g$ a
left-covariant $\sigma $-metric of the pair $(\Gp ,\Gm )$.
Then there are precisely two possibilities:\\
(i) Both $\Gp $ and $\Gm $ contain a unique (up to a constant factor)
non-zero left-invariant form of (the same) maximal degree.\\
(ii) Both $\Gp $ and $\Gm $ are infinite dimensional and for any
form $\omega \in \linv{\Gaw[k]{\tau }}$ there is a one-form
$\rho \in \Ga{\tau }$ such that $\omega \wedge \rho \not=0$.
\end{folg}

Let us fix $\xi '_{n_0}=\omega _0^{\tau _0}$ and
$\xi _{n_0}\in \Gaw[n_0]{-\tau _0}$ such that
$\ctrp{\omega _0^{\tau _0}}{\xi _{n_0}}=1$. {}From
Proposition \ref{s-hfumfr} we obtain
$\ctrpm{\omega _0^{\tau _0}}{\xi _{n_0}}\omega _0^{\tau _0}=
\omega _0^{\tau _0}\ctrmp{\xi _{n_0}}{\omega _0^{\tau _0}}$.
Hence the numbers
$\ctrpm{\omega _0^{\tau _0}}{\xi _{n_0}}$ and
$\ctrmp{\xi _{n_0}}{\omega _0^{\tau _0}}$ coincide.
Since $\dim \linv{\Gaw[n_0]{-\tau _0}}=1$ by Theorem \ref{t-eindhf},
$\xi _{n_0}$ is an eigenvector of $\sigma _{(n_0)}$.
Let $\sigma _{(n_0)}\xi _{n_0}=\lambda \xi _{n_0}$. Then we conclude
from the definition of the contractions and the considerations above that
\begin{align*}
1&=\ctrp{\omega _0^{\tau _0}}{\xi _{n_0}}
=\ctrm{\omega _0^{\tau _0}}{(-1)^{n_0(n_0-1)/2}\sigma _{(n_0)}\xi _{n_0}}\\
&=(-1)^{n_0(n_0-1)/2}\lambda \ctrm{\omega _0^{\tau _0}}{\xi _{n_0}}
=(-1)^{n_0(n_0-1)/2}\lambda \ctrp{\xi _{n_0}}{\omega _0^{\tau _0}}\\
&=(-1)^{n_0(n_0-1)/2}\lambda
\ctrm{(-1)^{n_0(n_0-1)/2}\sigma _{(n_0)}\xi _{n_0}}{\omega _0^{\tau _0}}
=\lambda ^2 \ctrm{\xi _{n_0}}{\omega _0^{\tau _0}}\\
&=\lambda ^2 \ctrp{\omega _0^{\tau _0}}{\xi _{n_0}}=\lambda ^2.
\end{align*}
Therefore, $\sigma _{(n_0)}^2(\xi _{n_0})=\xi _{n_0}$.

A consequence of Theorem \ref{t-eindhf} is that
$\kowr (\omega _0^{\tau _0})=\omega _0^{\tau _0}\otimes v_0$,
where $v_0$ is a one-dimensional corepresentation of $\A $. By assumption
(I) stated at the beginning of this section, it follows that
$v_0=1$. Hence $\omega _0^{\tau _0}$ is biinvariant.
Similarly, $\xi _{n_0}$ is biinvariant.
This implies that
$\sigma ^-_{(n_0,n_0)}(\omega _0^{\tau _0}\otA \xi _{n_0})=
\xi _{n_0}\otA \omega _0^{\tau _0}$. Applying Lemma \ref{l-symmet} we get
\begin{align*}
\ctrp{\xi _{n_0}}{\omega _0^{\tau _0}}&=
\ctrp{\cdot }{\cdot }\bigl( \sigma ^-_{(n_0,n_0)}(\omega _0^{\tau _0}\mt
\xi _{n_0})\bigr) \\
&=\ctrp{\omega _0^{\tau _0}}{\sigma ^2_{(n_0)}(\xi _{n_0})}
=\ctrp{\omega _0^{\tau _0}}{\xi _{n_0}}=1.
\end{align*}
This means that $\ctrp{\xi _{n_0}}{\omega _0^{\tau _0}}=1$ and hence
$\ctrm{\omega _0^{\tau _0}}{\xi _{n_0}}=
\ctrp{\xi _{n_0}}{\omega _0^{\tau _0}}=1$ and
$\ctrm{\xi _{n_0}}{\omega _0^{\tau _0}}=
\ctrp{\omega _0^{\tau _0}}{\xi _{n_0}}=1$.

Further, we have $\xi _{n_0}\wedge \rho =0$ for all $\rho \in \Ga{-\tau _0}$.
Therefore, the triple $(n_0,-\tau _0,\xi _{n_0})$ satisfies assumption (II)
at the beginning of the section as well. Now we can set
$\omega _0^{-\tau _0}:=\xi _{n_0}$ and so (\ref{eq-ctrnorm}) is valid.
In particular, we have obtained that
\begin{gather}\label{eq-sighf}
\sigma _{(n_0)}\omega _0^\pm =(-1)^{n_0(n_0-1)/2}\omega _0^\pm .
\end{gather}

Since the triple $(n_0,-\tau _0,\omega _0^{-\tau _0})$ satisfies
assumption (II), we can replace $\tau _0$ by $-\tau _0$
and Proposition \ref{s-hfumfr} remains true.
Moreover, it follows from Theorem
\ref{t-eindhf} that $\rho \wedge \omega _0^\pm =0$ for all
$\rho \in \Gaw[n_0]{\pm }$. Using this ansatz a similar reasoning as
used in the proof of Proposition \ref{s-hfumfr} shows the following.

\begin{satz}\label{s-hfumfl}
For any $\xi _k\in \Gaw[k]{-\tau }$ and $\xi '_l\in \Gaw[l]{\tau }$,
$0\leq l\leq k\leq n_0$, $\tau \in \{+,-\}$ the equations
\begin{align}
\xi '_l\wedge \ctrpm{\xi _k}{\omega _0^\tau }=
\ctrpm{\ctrmp{\xi '_l}{\xi _k}}{\omega _0^\tau }
\end{align}
hold.
\end{satz}

Let $\HoLpm ,\HoRpm :\Gaw[k]{\tau }\to \Gaw[n_0-k]{-\tau }$ denote the maps
given by
\begin{align}\label{eq-Hodef}
\HoLpm (\xi):=\ctrpm{\xi }{\omega _0^{-\tau }},\quad
\HoRpm (\xi):=\ctrpm{\omega _0^{-\tau }}{\xi }
\end{align}
for any $\xi \in \Gaw[k]{\tau }$, $0\leq k\leq n_0$, $\tau \in \{+,-\}$.

\begin{lemma}\label{l-Hoprops}
(i) For any $a\in \A $ and $\xi \in \Gaw{\tau }$, $\tau \in \{+,-\}$ we have
$\HoLpm (a\xi )=a\HoLpm (\xi )$ and $\HoRpm (\xi a)=\HoRpm (\xi )a$.\\
(ii) $\HoLp \HoLm =\HoLm \HoLp =\id $ and
$\HoRp \HoRm =\HoRm \HoRp =\id $.
In particular, the mappings $\HoLpm $ and $\HoRpm $ are isomorphisms of
$\Gaw{\tau }$ and $\Gaw{-\tau }$ as left and right $\A $-modules,
respectively.\\
(iii) For any $\rho _i\in \Gaw[k_i]{\tau }$, $i=1,2$, $k_1+k_2\leq n_0$,
$\tau \in \{+,-\}$, we have
\begin{gather}\label{eq-HoLRw}
\HoLpm (\rho _1\wedge \rho _2)=\ctrpm{\rho _1}{\HoLpm (\rho _2)},\quad
\HoRpm (\rho _1\wedge \rho _2)=\ctrpm{\HoRpm (\rho _1)}{\rho _2},\\
\label{eq-HorHol}
\ctrpm{\rho _1}{\HoRpm (\rho _2)}=\ctrpm{\HoLpm (\rho _1)}{\rho _2}.
\end{gather}
\end{lemma}

\begin{bew}
Since $\ctrpm{\cdot }{\cdot }$ is a homomorphism of $\A $-bimodules,
(i) follows from (\ref{eq-Hodef}). (ii) is obtained from
Proposition \ref{s-hfumfr} by inserting $\xi _k=\omega ^{-\tau }_0$ and
applying (\ref{eq-ctrnorm}).
Setting $\xi _0=\omega ^{-\tau }_0$ in Lemma \ref{l-ctrrek},
(\ref{eq-HoLRw}) and (\ref{eq-HorHol}) are equivalent to the equations
of Lemma \ref{l-ctrrek}(i) and \ref{l-ctrrek}(ii), respectively.
\end{bew}

\begin{defin}
We call the mapping $\HoLp :\Gaw{\tau }\to \Gaw{-\tau }$
\textit{left Hodge operator}
and $\HoRp :\Gaw{\tau }\to \Gaw{-\tau }$
\textit{right Hodge operator} on $\Gaw{\tau }$, $\tau \in \{+,-\}$.
\end{defin}

\begin{bem}
The equations in Proposition \ref{s-hfumfr} and \ref{s-hfumfl} with $k=l$
can also be written in the familiar form
\begin{align}
\HoRpm (\xi _k)\wedge \xi '_k&=\omega ^\tau _0\ctrmp{\xi _k}{\xi '_k},\\
\xi '_k\wedge \HoLpm (\xi _k)&=\ctrmp{\xi '_k}{\xi _k}\omega ^\tau _0.
\end{align}
\end{bem}

Up to now $\Gaw{+}$ and $\Gaw{-}$ have been only the exterior algebras
over bicovariant $\A $-bimodules $\Gp $ and $\Gm $, respectively.
In the remainder of this paper we assume in addition that
$\Gaw{\tau }$ is an inner bicovariant differential calculus with
differentiation $\dif _\tau $, $\tau \in \{+,-\}$. That the differential
calculus $\Gaw{\tau }$ is inner means that there exists a biinvariant
1-form $\eta ^\tau \in \Ga{\tau }$ such that
\begin{align}\label{eq-diff}
\dif _\tau \rho =\eta ^\tau \wedge \rho - (-1)^k\rho \wedge \eta ^\tau
\quad \rho \in \Gaw[k]{\tau },\tau \in \{+,-\}.
\end{align}
Further, we assume that the corresponding
$\sigma $-metrics (and hence contractions) are left-covariant.

\begin{defin}\label{d-kodiff0}
The mappings $\kodLpm :\Gaw[k]{\tau }\to \Gaw[k-1]{\tau }$
defined by
\begin{gather*}
\kodLpm \rho :=(-1)^k\HoLpm (\dif _{-\tau }\HoLmp (\rho )),\quad
\rho \in \Gaw[k]{\tau }, 0\leq k\leq n_0, \tau \in \{+,-\}
\end{gather*}
are called \textit{(positive and negative) left codifferential
operators on} $\Gaw{\tau }$.
Analogously we define the \textit{right codifferential operators}
$\kodRpm :\Gaw[k]{\tau }\to \Gaw[k-1]{\tau }$, $0\leq k\leq n_0$,
$\tau \in \{+,-\}$ \textit{on} $\Gaw{\tau }$ by
$\kodRpm \rho :=(-1)^{n_0-1+k}\HoRpm (\dif _{-\tau }\HoRmp (\rho ))$.
\end{defin}

\begin{lemma}\label{l-Hopmgleich}
$\HoLp (\rho )=\HoLm (\rho )$ and
$\HoRp (\rho )=\HoRm (\rho )$ for any $\rho \in \Gaw[k]{\tau }$,
$\tau \in \{+,-\}$, $k\in \{0,1,n_0-1,n_0\}$.
\end{lemma}

\begin{bew}
For $k=0$ we have $\HoLp (\rho )=\HoLm (\rho )=\rho \omega ^{-\tau }_0$
and $\HoRp (\rho )=\HoRm (\rho )=\omega ^{-\tau }_0\rho $ by definition.
For $k=n_0$ we obtain from Theorem \ref{t-eindhf} that there are $a,b\in \A $
such that $\rho =a\omega ^\tau _0=\omega ^\tau _0b$.
Then Lemma \ref{l-Hoprops}(i) and equation (\ref{eq-ctrnorm}) imply that
$\HoLpm (a\omega ^\tau _0)=a\HoLpm (\omega ^\tau _0)=
a\ctrpm{\omega ^\tau _0}{\omega ^{-\tau }_0}=a$ and
$\HoRpm (\omega ^\tau _0b)=\HoLpm (\omega ^\tau _0)b=
\ctrpm{\omega ^{-\tau }_0}{\omega ^\tau _0}b=b$.

Let now $k=n_0-1$. We compute
\begin{gather*}\tag{$*$}
\begin{aligned}
\HoLpm (\rho )&=\ctrpm{\rho }{\omega ^{-\tau }_0}
=\tilde{g}(A^\pm _{n_0-1}\rho \mt B^\pm _{n_0-1,1}\omega ^{-\tau }_0)\\
&=\tilde{g}\bigl( \rho \mt (A^\pm _{n_0-1}\otA \id)B^\pm _{n_0-1,1}
\omega ^{-\tau }_0\bigr)
=\tilde{g}(\rho \mt A^\pm _{n_0}\omega ^{-\tau }_0)
\end{aligned}
\end{gather*}
by using Lemma \ref{l-tgant} and the second equation of (\ref{eq-AB}).
We also have $A^+_k=(-1)^{k(k-1)/2}A^-_k\sigma ^+_{(k)}$ for any $k\geq 1$.
Hence (\ref{eq-sighf}) gives
\begin{align*}
A^+_{n_0}\omega ^{-\tau }_0&=(-1)^{n_0(n_0-1)/2}A^-_{n_0}\sigma ^+_{(n_0)}
\omega ^{-\tau }_0\\
&=(-1)^{n_0(n_0-1)/2}A^-_{n_0}(-1)^{n_0(n_0-1)/2}\omega ^{-\tau }_0
=A^-_{n_0}\omega ^{-\tau }_0.
\end{align*}
{}From this and equation ($*$) we conclude that $\HoLp (\rho )=\HoLm (\rho )$.

In the case $k=1$ we use that the mappings $\HoLpm $ are isomorphisms of left
$\A $-modules. Therefore there is a $\rho '\in \Gaw[n_0-1]{\tau }$ such that
$\rho = \HoLp (\rho ')$. By the preceding we also have $\rho =\HoLm (\rho ')$.
Hence, $\HoLp (\rho )=\HoLp \HoLm (\rho ')=\rho '$ and
$\HoLm (\rho )=\HoLm \HoLp (\rho ')=\rho '$. Similarly,
$\HoRp (\rho )=\HoRm (\rho )$ for any $\rho \in \Gaw[k]{\tau }$,
$k=1,n_0-1$.
\end{bew}

\begin{lemma}\label{l-HoLeqR}
For any $\rho \in \rinv{\Ga{\tau}}$, $\tau \in \{+,-\}$
we have $\HoLpm (\rho )=(-1)^{n_0-1}\HoRpm (\rho )$.
\end{lemma}

\begin{bew}
The $n_0$-form $\omega _0^{-\tau }$ is left-invariant. Hence there are
left-invariant 1-forms $\rho _1,\ldots,\rho _{n_0}\in \linv{\Ga{-\tau}}$
such that $\omega _0^{-\tau}=\rho _1\wedge \ldots \wedge \rho _{n_0}$.
Then (\ref{eq-ctrrekl}) and the $\sigma $-symmetry of the $\sigma $-metric
yield
\begin{align}
\HoLp (\rho )=\sum _{i=1}^{n_0}(-1)^{i-1}\rho _1\wedge \ldots \wedge \rho _{i-1}
\ctr{\rho }{\rho _i}\wedge \rho _{i+1}\wedge \ldots \wedge \rho _{n_0}.
\end{align}
The $\sigma $-symmetry of the $\sigma $-metric implies that
$\ctr{\rho _i}{\rho }=\ctr{\rho }{\rho _i}$ for any $i=1,\ldots ,n_0$.
Using this fact and equation (\ref{eq-ctrrekr}) we obtain
the same formula for $(-1)^{n_0-1}\HoRm (\rho )$.
Applying Lemma \ref{l-Hopmgleich} the assertion follows.
\end{bew}

\begin{satz}\label{s-kodLR}
The codifferentials $\kodL{\tau '}$ and $\kodR{\tau '}$,
$\tau '\in \{+,-\}$, coincide.
On $a\in\A$ they act trivially: $\kodLpm a=\kodRpm a=0$. For
any $\rho\in\Gaw[k]{\tau}$, $k>0$, $\tau \in\{+,-\}$ we have
\begin{gather}\label{eq-kodpm0}
\kodLpm \rho =\ctrpm{\rho }{\eta ^{-\tau }}
 +(-1)^k\ctrpm{\eta ^{-\tau }}{\rho }.
\end{gather}
\end{satz}

\begin{bew}
Let $k>0$ and $\rho \in \Gaw[k]{\tau }$. The definition of
$\kodLpm$ and (\ref{eq-diff}) give
\begin{align*}
\kodLpm \rho &=(-1)^k\HoLpm (\dif _{-\tau } \HoLmp (\rho ))=
(-1)^k\HoLpm (\eta ^{-\tau } \wedge \HoLmp (\rho )
-(-1)^{n_0-k}\HoLmp (\rho )\wedge \eta ^{-\tau }).
\end{align*}
{}From the first equation of (\ref{eq-HoLRw}) and Lemma \ref{l-Hoprops}(ii)
we obtain that the first summand is equal to
$(-1)^k\ctrpm{\eta ^{-\tau }}{\HoLpm (\HoLmp (\rho ))}=
(-1)^k\ctrpm{\eta ^{-\tau }}{\rho }$.
For the second summand we use (\ref{eq-HoLRw}) and Lemma \ref{l-HoLeqR}
and obtain
$(-1)^{n_0+1}\ctrpm{\HoLmp (\rho )}{\HoLpm (\eta ^{-\tau })}=
\ctrpm{\HoLmp (\rho )}{\HoRpm (\eta ^{-\tau })}$.
We apply now (\ref{eq-HorHol}) and Lemma \ref{l-Hoprops}(ii) to the latter
and get
$\ctrpm{\HoLpm (\HoLmp (\rho ))}{\eta ^{-\tau }}=
\ctrpm{\rho }{\eta ^{-\tau }}$. This proves (\ref{eq-kodpm0}) for the left
codifferentials. Similar computations lead to the same expression
for $\kodRpm \rho $.
\end{bew}

\begin{satz}
For any $\rho \in \linv{\Gaw[n_0-1]{\tau }}$, $\tau \in \{+,-\}$ we have
$\dif _\tau \rho =0$.
\end{satz}

\begin{bew}
Let $\rho \in \linv{\Gaw[n_0-1]{\tau }}$. Because of Lemma
\ref{l-Hoprops}(ii) and the left-covariance of $\HoLpm $ there are
$\rho ^\pm _1\in \linv{\Ga{-\tau }}$ such that
$\rho =\HoLpm (\rho ^\pm _1)$. Then $\dif _\tau \rho =0$
is equivalent to
\[
0=\HoLpm (\dif _\tau \rho )=\HoLpm (\dif _\tau \HoLmp (\rho ^\mp _1))
=-\kodLpm \rho ^\mp _1.
\]
Since $\eta ^\tau $ is biinvariant, $\rho ^\mp _1$ is left-invariant
and the $\sigma $-metric is $\sigma $-symmetric,
we conclude from Proposition \ref{s-kodLR} that
\[
\kodLpm \rho ^\mp _1=\ctrpm{\rho ^\mp _1}{\eta ^\tau }
-\ctrpm{\eta ^\tau }{\rho ^\mp _1}
=\ctrpm{\rho ^\mp _1}{\eta ^\tau }-\ctrpm{\rho ^\mp _1}{\eta ^\tau }=0.
\]
\end{bew}

\section{Laplace-Beltrami Operators}
\label{sec-Laplace}

Let $\A$ be again an arbitrary Hopf algebra and let $\Gp ,\Gm $ be
two bicovariant $\A $-bimodules which admit a left-covariant $\sigma $-metric
in the sense of Definition \ref{d-metrik}. Moreover, (as in the
last part of Section \ref{sec-Hodge},) we assume that
the bicovariant $\A $-bimodules $\Gaw{\tau }$, $\tau \in \{+,-\}$
admit a differential operator $\dif _\tau $ such that they become inner
bicovariant differential calculi on $\A $.
Further we suppose that the $\sigma $-metrics (and hence contractions)
are left-covariant.

In addition we now assume that the Hopf algebra $\A $ is cosemisimple
\cite[Sect.\ 11.2]{b-KS}, that is, there exists a linear functional $h$
on $\A $, called the Haar functional, such that $h(1)=1$ and
\begin{gather}\label{eq-Haar}
(h\otimes \id )\komult (a)=(\id \otimes h)\komult (a)=h(a)1
\end{gather}
for all $a\in \A $. Further, we suppose that the Haar functional
is regular, that is,
both $h(ab)=0$ for all $b\in \A $ and $h(ba)=0$ for all $b\in \A $
imply that $a=0$. (Recall that any CQG-algebra is cosemisimple and its
Haar functional is regular (\cite{b-KS}, Proposition 11.29).
By Proposition \ref{s-metna} the restriction of
$\ctrpm{\cdot }{\cdot }$ to $\Gaw[k]{\tau }\tet \Gaw[k]{-\tau }$ is
nondegenerate. Hence for each $\rho \in \Gaw[k]{\tau }$ there is a
$\rho '\in \Gaw[k]{-\tau }$ such that
$\A \ni a:=\ctrpm{\rho }{\rho '}\not=0$.
By the regularity of the Haar functional there is a $b\in \A $ such that
$h(ab)\not=0$. Then we have $h\ctrpm{\rho }{\rho 'b}=
h(\ctrpm{\rho }{\rho '}b)=h(ab)\not=0$. Therefore, the mapping
$h\circ \ctrpm{\cdot }{\cdot }:\Gaw[k]{\tau }\tet \Gaw[k]{-\tau }\to \comp $
is nondegenerate for all $k\geq 0$ and $\tau \in \{+,-\}$.
We shall consider it as a generalisation of the classical notion of the
metric on $k$-forms.

Motivated by Definition \ref{d-kodiff0} and Proposition \ref{s-kodLR},
we introduce the following notion.

\begin{defin}\label{d-kodiff}
The mappings $\kodpm{\tau } :\Gaw[k]{\tau }\to \Gaw[k-1]{\tau }$, $k\geq 0$,
$\tau \in \{+,-\}$, defined by $\kodpm{\tau }(a)=0$ for $a\in \A $ and
\begin{gather}\label{eq-kodpm}
\kodpm{\tau }\rho =\ctrpm{\rho }{\eta ^{-\tau }}
 +(-1)^k\ctrpm{\eta ^{-\tau }}{\rho }
\end{gather}
for $\rho \in \Gaw[k]{\tau }$, $k>0$, are called
\textit{(positive and negative) codifferential operators on} $\Gaw[k]{\tau }$.
\end{defin}

\begin{lemma}\label{l-kodeig}
(i) $(\kodpm{\tau })^2=0$.\\
(ii) $\kodpm{\tau }(a\rho )=a\kodpm{\tau }\rho +
(-1)^k\ctrpm{\dif_{-\tau }a}{\rho }$ for any $a\in \A $,
$\rho \in \Gaw[k]{\tau }$, $\tau \in \{+,-\},k\geq 1$.
\end{lemma}

\begin{bew}
(i) Since $(\kodpm{\tau })^2(\rho )\in \Gaw[k-2]{\tau }$ for any
$\rho \in \Gaw[k]{\tau }$, $k\geq 0$, $\tau \in \{+,-\}$, we obtain
$(\kodpm{\tau })^2(\rho )=0$ for $\rho \in \Gaw[k]{\tau }$, $k\leq 1$.
For $k\geq 2$ we get
\begin{align*}
(\kodpm{\tau })^2(\rho )=&
\kodpm{\tau }\bigl( \ctrpm{\rho }{\eta ^{-\tau }}
+(-1)^k\ctrpm{\eta ^{-\tau }}{\rho }\bigr) \\
=&
\ctrpm{
\bigl( \ctrpm{\rho }{\eta ^{-\tau }}+(-1)^k\ctrpm{\eta ^{-\tau }}{\rho }\bigr)
}{\eta ^{-\tau }}\\
&+(-1)^{k-1}\ctrpm{\eta ^{-\tau }}{
\bigl( \ctrpm{\rho }{\eta ^{-\tau }}+(-1)^k\ctrpm{\eta ^{-\tau }}{\rho }\bigr)
}.\\
\intertext{Applying Lemma \ref{l-ctrrek}(i) on the first and fourth summand
we obtain}
=&
\ctrpm{\rho }{\eta ^{-\tau }\wedge \eta ^{-\tau }}
+(-1)^k\ctrpm{\ctrpm{\eta ^{-\tau }}{\rho }}{\eta ^{-\tau }}\\
&+(-1)^{k-1}\ctrpm{\eta ^{-\tau }}{\ctrpm{\rho }{\eta ^{-\tau }}}
-\ctrpm{\eta ^{-\tau }\wedge \eta ^{-\tau }}{\rho }.
\end{align*}
Since $\eta ^{-\tau }$ is biinvariant, $\eta ^{-\tau}\wedge \eta ^{-\tau }=0$.
Using Lemma \ref{l-ctrrek}(ii) the second and third summand in the last
expression also vanish.

(ii) From (\ref{eq-diff}) it follows that
\[
\ctrpm{\dif _{-\tau }a}{\rho }=\ctrpm{\eta ^{-\tau }a}{\rho }
-\ctrpm{a\eta ^{-\tau }}{\rho }
=\ctrpm{\eta ^{-\tau }}{a\rho }-a\ctrpm{\eta ^{-\tau }}{\rho }.
\]
Then (\ref{eq-kodpm}) gives the assertion.
\end{bew}

\begin{lemma}\label{l-hmetprop}
For any $a\in \A $ and $\rho \in \linv{\Ga{\tau }},
\rho '\in \linv{\Ga{-\tau }}$, $\tau \in \{+,-\}$ we have\\
(i) $h(\ctrpm{a\rho }{\rho '})=h(\ctrpm{\rho a}{\rho '})
=h(a)\ctrpm{\rho }{\rho '}$,\\
(ii) $h\bigl( \kodpm{-\tau }(a\rho ')\bigr)=0$.
\end{lemma}

\begin{bew}
(i) Let $\{\theta _i\,|\,i=1,\ldots,m\}$ be a basis of the vector space
$\linv{\Ga{\tau}}$. It suffices to prove the assertion for $\rho =\theta _i$.
The left-invariance of the $\sigma $-metric ensures that
$\ctrpm{\rho }{\rho '}\in \comp$ and we conclude that
$h(\ctrpm{a\rho }{\rho '})=h(a\ctrpm{\rho }{\rho '})
=h(a)\ctrpm{\rho }{\rho '}$.

By the general theory \cite{a-Woro2} there are functionals $f^i_j$,
$i,j=1,\ldots,m$,
such that $\theta _ia=a_{(1)}f^i_j(a_{(2)})\theta _j$ and
$f^i_j(1)=\delta^i_j$. We have again
$\ctrpm{\theta _j}{\rho '}\in \comp $ and therefore
\begin{align*}
h(\ctrpm{\theta _ia}{\rho '})
&=h(a_{(1)})f^i_j(a_{(2)})\ctrpm{\theta _j}{\rho '}
=f^i_j(h(a_{(1)})a_{(2)})\ctrpm{\theta _j}{\rho '}\\
&=f^i_j(h(a)\cdot 1)\ctrpm{\theta _j}{\rho '}
=h(a)\ctrpm{\theta _i}{\rho '}
\end{align*}
by (\ref{eq-Haar}). Hence we get (i).

(ii) Firstly we see from (\ref{eq-kodpm}) that
$\kodpm{-\tau }(\rho ')=\ctrpm{\rho '}{\eta ^\tau }-
\ctrpm{\eta ^\tau }{\rho '}=0$ since the $\sigma $-metric is
$\sigma $-symmetric,
$\eta ^\tau $ is biinvariant and $\rho '$ is left-invariant.
Secondly, Lemma \ref{l-kodeig}(ii) gives
$h(\kodpm{-\tau }(a\rho '))=h(a\kodpm{-\tau }\rho '
-\ctrpm{\dif _\tau a}{\rho '})
=h\bigl( \ctrpm{a\eta ^\tau }{\rho '}-\ctrpm{\eta ^\tau a}{\rho '}\bigr)$.
Then the assertion follows from (i).
\end{bew}

\begin{thm}\label{t-hdifkod}
Suppose that $g$ is a left-invariant $\sigma $-metric of the pair
$(\Gp ,\Gm )$. Let
$\ctrpm{\cdot }{\cdot }$ be the corresponding contractions.
Then for any $\rho \in \Gaw[k]{\tau }$, $\rho '\in \Gaw[k+1]{-\tau }$,
$\tau \in \{+,-\}$ the equations
\begin{align}\label{eq-hdifkod}
h(\ctrpm{\rho }{\kodpm{-\tau }\rho '})&=h(\ctrpm{\dif _\tau \rho }{\rho '})
\quad \text{and}\\
\label{eq-hkoddif}
h(\ctrpm{\kodpm{-\tau }\rho '}{\rho })&=h(\ctrpm{\rho '}{\dif _\tau \rho })
\end{align}
hold.
\end{thm}

\begin{bew}
Inserting the definitions (\ref{eq-kodpm}) and (\ref{eq-diff}) we obtain
\begin{align*}
&h\bigl(
\ctrpm{\rho }{\kodpm{-\tau }\rho '}-\ctrpm{\dif _\tau \rho }{\rho '}\bigr)
=h\Bigl(
\ctrpm{\rho }{\bigl(
\ctrpm{\rho '}{\eta ^\tau }+(-1)^{k+1}\ctrpm{\eta ^\tau }{\rho '}\bigr)}\\
&\quad -\ctrpm{(\eta ^\tau \wedge \rho
+(-1)^{k+1}\rho \wedge \eta ^\tau )}{\rho '} \Bigr).
\end{align*}
Applying Lemma \ref{l-ctrrek} we now substitute
$\ctrpm{\rho }{\ctrpm{\rho '}{\eta ^\tau }}$ by
 $\ctrpm{\ctrpm{\rho }{\rho '}}{\eta ^\tau }$,
$\ctrpm{\rho }{\ctrpm{\eta ^\tau }{\rho '}}$ by
 $\ctrpm{\rho \wedge \eta ^\tau }{\rho '}$ and
$\ctrpm{\eta ^\tau \wedge \rho }{\rho '}$ by
 $\ctrpm{\eta ^\tau }{\ctrpm{\rho }{\rho '}}$.
Then we have
\begin{align*}
h(\ctrpm{\rho}{\kodpm{-\tau}\rho '})-
h(\ctrpm{\dif _\tau \rho }{\rho '})&=
h\bigl( \ctrpm{\ctrpm{\rho }{\rho '}}{\eta ^\tau}
-\ctrpm{\eta ^\tau }{\ctrpm{\rho }{\rho '}}\bigr)\\
&=h\bigl( \kodpm{-\tau }\ctrpm{\rho }{\rho '}\bigr).
\end{align*}
Since $\ctrpm{\rho }{\rho '}$ is an element of
$\Ga{-\tau }=\A \linv{\Ga{-\tau }}$, we obtain (\ref{eq-hdifkod})
by Lemma \ref{l-hmetprop}(ii). The proof of (\ref{eq-hkoddif}) is similar.
\end{bew}

\begin{defin}\label{d-Laplace}
We call the operators $\Lappm :\Gaw[k]{\tau }\to \Gaw[k]{\tau }$,
$\Lappm :=\dif _\tau \kodpm{\tau }+ \kodpm{\tau }\dif _\tau $
\textit{Laplace-Beltrami operators}.
\end{defin}

The following properties of $\Lappm $ are simple consequences of
the facts that $\dif ^2=0$, $(\kodpm{\tau })^2=0$ and (\ref{eq-kodpm}).

\begin{lemma}\label{l-Lapprop}
The Laplace-Beltrami operators satisfy the equations
\begin{align}
\label{eq-LapDir}
\Lappm &=(\dif _\tau +\kodpm{\tau })^2,\\
\label{eq-Lapdif}
\Lappm \dif _\tau &=\dif _\tau \Lappm =\dif _\tau \kodpm{\tau }\dif _\tau ,\\
\label{eq-Lapkod}
\Lappm \kodpm{\tau }&=\kodpm{\tau }\Lappm
=\kodpm{\tau }\dif _\tau \kodpm{\tau },\\
\label{eq-LapA}
\Lap[+]{\tau '}a=\Lap[-]{\tau '}a&=\ctr{\eta^+a}{\eta^-}
+\ctr{\eta^-a}{\eta^+}-2a\ctr{\eta^+}{\eta^-}
\end{align}
for any $a\in \A $ and $\tau ,\tau '\in\{+,-\}$.
\end{lemma}

\begin{bem}
By (\ref{eq-LapA}) the Laplace-Beltrami operator on $\A \subset \Gaw{\pm }$
neither depends on the sign $\tau '$ of the antisymmetrizer
nor on the $\A $-bimodule $\Gaw{\pm }$ containing $\A $.
\end{bem}

\begin{satz}
For any $\rho \in \Gaw[k]{\tau }$, $\rho '\in \Gaw[k]{-\tau }$,
$\tau \in \{+,-\}$, $k\geq 0$ we have
\begin{gather}
h(\ctrpm{\Lappm \rho }{\rho '})=h(\ctrpm{\rho }{\Lappm[-\tau ]\rho '}).
\end{gather}
\end{satz}

\begin{bew}
Using Theorem \ref{t-hdifkod} we compute
\begin{align*}
h(\ctrpm{\Lappm \rho }{\rho '})&
=h(\ctrpm{\dif _\tau \kodpm{\tau }\rho }{\rho '}
+\ctrpm{\kodpm{\tau }\dif _\tau \rho }{\rho '})\\
&=h(\ctrpm{\kodpm{\tau }\rho }{\kodpm{-\tau }\rho '}
+\ctrpm{\dif _\tau \rho }{\dif _{-\tau }\rho '})\\
&=h(\ctrpm{\rho }{\dif _{-\tau }\kodpm{-\tau }\rho '}
+\ctrpm{\rho }{\kodpm{-\tau }\dif _{-\tau }\rho '})
=h(\ctrpm{\rho }{\Lappm[-\tau ]\rho '}).
\end{align*}
\end{bew}

\section{Eigenvalues of the Laplace-Beltrami operator for $\SLqN$}
\label{sec-eigenvalues}

Throughout this section we assume that $q$ is a transcendental complex
number and
$\A $ is the Hopf algebra $\OSLqN $, $N\geq 2$.
Then $\A $ is cosemisimple, i.\,e.\ any element of $\A $ is a finite
linear combination of matrix elements of irreducible matrix
corepresentations of $\A $ (\cite{b-KS}, Theorem 11.22).
Further, $\A $ is coquasitriangular and admits a universal $r$-form
$\br :\A \otimes \A \to \comp$ defined by
$\br (u^i_j\otimes u^k_l)=z^{-1}\Rda{}^{ki}_{jl}$, where $z$ is a fixed
complex number with $z^N=q$, and
\begin{gather}
\Rda{}^{ij}_{kl}=q^{\delta ^i_j}\delta ^i_l\delta ^j_k
+(i<j)(q-q^{-1})\delta ^i_k\delta ^j_l.
\end{gather}
Here the number $(i<j)$ is 1 if $i<j$ and zero otherwise. We shall write
$\Rda {}^\pm $ for $\Rda {}^{\pm 1}$.

Let $\Gp $ and $\Gm $ be the
$N^2$-dimensional bicovariant differential calculi on $\A $ determined
by the fundamental corepresentation $u$ and the contragredient
corepresentation $u\cont $ (see Section \ref{sec-examples}).
Further, let denote $F_1,F_2,G_1,G_2$ the $N\times N$-matrices with entries
$F_1{}^i_j=z^{-1}q^{N-2i}\delta ^i_j$, $F_2{}^i_j=q^{2i}\delta ^i_j$,
$G_1{}^i_j=z^{-1}q^N\delta ^i_j$, $G_2{}^i_j=\delta ^i_j$.
Then $F_1\in \Mor(u^{\contr\contr},u)$, $F_2\in \Mor(u,u^{\contr\contr})$
and $G_1,G_2\in \Mor(u)$ and they determine a bicovariant $\sigma $-metric
of the pair $(\Gp,\Gm)$ (see Section \ref{sec-examples}).

The Laplace-Beltrami operator $\Lap[]{}$ on $\A $ is given by
(\ref{eq-LapA}). For $n\in \mathbb{Z}$ and a complex number
$p\not=0,\pm 1$ let $[n]_p$ denote the number $(p^n-p^{-n})/(p-p^{-1})$.

\begin{satz}\label{s-Lapew}
The Laplace-Beltrami operator $\Lap[]{}$ on $\A $ is diagonalizable.
Let $v^\lambda $ be a fixed irreducible corepresentation of $\A $
corresponding to a Young diagram $\lambda $. Then the matrix elements
of $v^\lambda $ are eigenvectors of $\Lap[]{}$ to the eigenvalue
\begin{equation}\label{eq-Lapew}
E_\lambda :=(z-z^{-1})^2\left([m]_z^2[N]_q
+[N]_z\sum _{(i,j)\in \lambda }[N^2-2m+2N(j-i)]_z\right),
\end{equation}
where $(i,j)\in \lambda $ means that there is a box in the
$i$-th row and $j$-th column of $\lambda $ and $m$ is the number
of boxes in $\lambda $.
\end{satz}

\begin{bew}
Using the relations $\br (u^i_j,S(u^k_l))=zq^{2k-2l}\Rdam {}^{ik}_{lj}$
and $\br (S(u^i_j),u^k_l)=z\Rdam {}^{ik}_{lj}$, some properties of
the $r$-form $\br $ and the $\Rda $-matrix equation (\ref{eq-LapA}),
for any $m\geq 0$ we get
\begin{align*}
\Lap[]{}(u^{i_1}_{j_1}u^{i_2}_{j_2}\ldots u^{i_m}_{j_m})&=
q^{-N-1}u^{i_1}_{l_1}u^{i_2}_{l_2}\ldots u^{i_m}_{l_m}\bigl(
z^{-2m}D^+_{m+1}+z^{2m}D^-_{m+1}-2\id \bigr)
^{l_1\ldots l_m k}_{j_1\ldots j_m n}q^{2k}\delta ^n_k,
\end{align*}
where
\begin{gather}\label{eq-Jucys}
D^\pm _{m+1}=\Rda {}^\pm _{m,m+1}\Rda {}^\pm _{m-1,m}\ldots
\Rda {}^\pm _{23}\Rda {}^\pm _{12}{}^2\Rda {}^\pm _{23}\ldots
\Rda {}^\pm _{m,m+1}, \quad m\geq 2
\end{gather}
are the so called Jucys-Murphy operators of the Hecke algebra, $D^\pm _1=\id $.
Since $q^{2d-2a}\Rda {}^{bd}_{ac}\Rda {}^{ec}_{fd}=\delta ^e_a\delta ^b_f
+q^{2N+1}(q-q^{-1})q^{-2a}\delta ^b_a\delta ^e_f$ and
$q^{2d-2a}\Rdam {}^{bd}_{ac}\Rdam {}^{ec}_{fd}=\delta ^e_a\delta ^b_f
-q(q-q^{-1})q^{-2a}\delta ^b_a\delta ^e_f$ we obtain
\begin{gather*}
\sum _k q^{2k}(D^\pm _{m+1})^{l_1l_2\ldots l_mk}_{j_1j_2\ldots j_mk}=
\bigl(q^{N+1}[N]_q\id \pm q^{N+1}q^{\pm N}(q-q^{-1})\sum _{n=1}^mD^\pm _n
\bigr)^{l_1\ldots l_m}_{j_1\ldots j_m}
\end{gather*}
and hence
\begin{align*}
\Lap[]{}(u^{i_1}_{j_1}\cdots u^{i_m}_{j_m})=u^{i_1}_{l_1}\cdots u^{i_m}_{l_m}
\bigl(&(z^m-z^{-m})^2[N]_q\id \\
&+(q-q^{-1})\sum _{n=1}^m(q^Nz^{-2m}D^+_n-
q^{-N}z^{2m}D^-_n)\bigr)^{l_1\ldots l_m}_{j_1\ldots j_m}.
\end{align*}
Since $q$ is transcendental, $\A $ is cosemisimple.
Moreover, $\A $ is generated by the matrix elements of the fundamental
corepresentation $u$ of $\A $.
Let $P_\lambda $ be a projection of $u^{\otimes m}$ onto the irreducible
corepresentation of $\A $ corresponding to the Young diagram $\lambda $.
Then Proposition 4.7 and the preceding considerations in \cite{a-LeducRam}
imply that
$\sum _{n=1}^mD^\pm _nP_\lambda =\sum _{(i,j)\in \lambda }q^{\pm (2j-2i)}
P_\lambda$ and therefore
$\Lap[]{}(u^{i_1}_{k_1}\cdots u^{i_m}_{k_m}
P_\lambda {}^{k_1\ldots k_m}_{j_1\ldots j_m})=E_\lambda
u^{i_1}_{k_1}\cdots u^{i_m}_{k_m}P_\lambda {}^{k_1\ldots k_m}_{j_1\ldots j_m}$,
where
\begin{gather*}
E_\lambda =(z^m-z^{-m})^2[N]_q+(q-q^{-1})\sum _{(i,j)\in \lambda }
(z^{-2m}q^{N+2j-2i}-z^{2m}q^{N-2j+2i}).
\end{gather*}
Since $q=z^N$, (\ref{eq-Lapew}) follows.
\end{bew}

\begin{bems}
1. The corepresentation $v^\lambda $ of $\A $ with Young diagram $\lambda $
corresponds to the representation of $U_q(\mf{g})$ with highest weight
$\lambda =\sum _{i=1}^{N-1}m_i\omega _i$ where $\omega _i$ are the
fundamental weights and $m_i$ the number of columns in $\lambda $ of
length $i$.
Let $\paar{\cdot }{\cdot }$ denote the Killing metric on the Lie
algebra $\mathrm{sl}_{N-1}$ and let $\rho _0$ be the half sum of positive
roots. Then the eigenvalue of the classical Laplace-Beltrami operator
(with respect to the biinvariant metric)
corresponding to the highest weight $\lambda $ is given by the formula
\begin{gather}\label{eq-cllap}
\tilde{E}_\lambda =\paar{\lambda +\rho _0}{\lambda +\rho _0}-
\paar{\rho _0}{\rho _0}=
\sum_{i=1}^{N-1}\frac{(N-i)m_i}{N}\left( i(m_i+N)
+2\sum_{j=1}^{i-1}jm_j\right )
\end{gather}
(see \cite{b-Wallach}).
For the quantum case one can check that
$\lim_{q\to 1}(q-1/q)^{-2}E_\lambda=\tilde{E}_\lambda$.

2. For $N=2$ we have $q=z^2$ and equation (\ref{eq-Lapew}) reduces
to the formula
\begin{equation}\label{eq-SU2ew}
E_{[m]} :=2(z-z^{-1})^2[m]_z[m+2]_z.
\end{equation}
\end{bems}

\begin{satz}
Let $z$ be a transcendental real number and $q=z^N$.\\
(i) All the eigenvalues of the
Laplace-Beltrami operator $\Lap[]{}:\A \to \A $ are nonnegative.\\
(ii) For any $a\in \A $ we have $\Lap[]{}(a)=0$ if and only if
$a\in \comp 1$.\\
(iii) The smallest positive eigenvalue of $\Lap[]{}:\A \to \A $ is
\begin{equation}\label{eq-kleinstew}
\mathrm{min}\{E_\lambda \,|\,\lambda =[1^k,0^{N-k}],k=1,\ldots,N-1\}.
\end{equation}
\end{satz}

\begin{bew}
We prove the assertions of the Proposition in the case $z>0$. The other cases
are an easy consequence of this one.

Firstly one shows that if $\lambda =[l_1,l_2,\ldots,l_N]$,
$l_1\geq l_2\geq \ldots \geq l_N\geq 1$, then $E_\lambda >E_{\lambda '}$,
where $\lambda '=[l_1-1,l_2-1,\ldots,l_N-1]$.
Secondly, if $\lambda =[l_1,l_2,\ldots ,l_k,0^{N-k}]$, $l_k>0$, $1\leq k<N$,
and $l_i>l_{i+1}$, $l_i\geq 2$ for some $i=1,2,\ldots ,k$, then let
$\lambda '$ be the diagram
$[l_1,l_2,\ldots ,l_{i-1},l_i-1,l_{i+1},\ldots ,l_k,1,0^{N-k-1}]$.
One can prove that $E_\lambda >E_{\lambda '}$ since $[n]_z>[n-2n']_z$
for all $n'\in \mathbb{N}$, $n\in \mathbb{Z}$. Therefore, for any
$\lambda \not=[0^N]$ there exists a $\lambda '=[1^k,0^{N-k}]$ such that
$E_\lambda \geq E_{\lambda '}$. Obviously, $E_{[0^N]}=0$ and
because of $[m]_p>0$ for any $m\in \mathbb{N}$, $p>0$, we also have 
\begin{align*}
E_{\lambda '}=(z-z^{-1})^2\bigl([k]_z^2[N]_q+
[N]_z[k]_q[(N+2)(N-k-1)+2]_z\bigr)>0
\end{align*}
for any $\lambda '=[1^k,0^{N-k}]$, $1\leq k<N$. Hence the assertions follow.
\end{bew}

\begin{bem}
Let $q$ be a transcendental complex number. Let $\A $ be one of the
quantum groups $\OSpqN $ or $\OOqN $, $N\geq 3$,
and $\Gp $, $\Gm $ as in Section \ref{sec-examples}, where $u$ is the
fundamental corepresentation of $\A $. Then the settings
$G_1{}^i_j:=\epsilon r/2\delta ^i_j$, $G_2{}^i_j:=\delta ^i_j$,
$F_1{}^i_j:=rq^{2\rho _i}/2\delta ^i_j$,
$F_2{}^i_j:=\epsilon q^{-2\rho _i}\delta ^i_j$,
where $r=\epsilon q^{N-\epsilon }$
(we use the notation of \cite{a-FadResTak1}),
determine a left-covariant $\sigma $-metric of the pair $(\Gp ,\Gm )$.
Similarly to the proof of Proposition \ref{s-Lapew}, using (6.14) in
\cite{a-LeducRam} one can show that the eigenvalues of the
Laplace-Beltrami operator $\Lap[]{}$ on $\A $ corresponding to the
Young diagram $\lambda $ are
\begin{gather*}
E_\lambda =(q-q^{-1})^2\sum _{(i,j)\in \lambda }[N-\epsilon +2j-2i]_q.
\end{gather*}
During the computations the operators
$r\sum _{k=1}^mD^+_k-r^{-1}\sum _{k=1}^mD^-_k$ of the Birman-Wenzl-Murakami
algebra appear ---
one can take (\ref{eq-Jucys}) for the definition of $D^\pm _k$, where
$\Rda {}^\pm $ denote the matrices
\begin{gather}
\Rda {}^\pm {}^{ij}_{kl}=q^{\pm (\delta ^i_j-\delta ^i_{j'})}
\delta ^i_l\delta ^j_k\pm (\pm i < \pm l)(q-q^{-1})(\delta ^i_k\delta ^j_l
-\epsilon _i\epsilon _lq^{\rho _l-\rho _i}\delta ^i_{j'}\delta ^k_{l'})
\end{gather}
---, which are central in the algebra $\Mor(u^{\otimes m+1})$.
\end{bem}

%\bibliography{quantum}
%\bibliographystyle{mybibger}

\end{document}